\documentclass[12pt]{article}
\usepackage{fullpage}
\usepackage{epsf}
\usepackage{epsfig}
\usepackage{amsthm}
\usepackage{amsfonts}
\usepackage{color}
\usepackage{amsmath}
\usepackage{setspace}
\usepackage{hyperref}
\usepackage{algorithmic}
\usepackage[boxed]{algorithm}
\usepackage[usenames, dvipsnames]{xcolor}

\usepackage{graphicx}
\usepackage{caption}
\usepackage{subcaption}
\usepackage[normalem]{ulem}
\usepackage{mathrsfs}

\usepackage{tikz}
\usetikzlibrary{shapes,arrows}

\newcommand{\blind}{0}

\begin{document}

\newcommand{\bm}[1]{\mbox{\boldmath $#1$}}
\newcommand{\mb}[1]{#1}%{\mathbf{#1}}
\newcommand{\E}{\mathbb{E}}
\newcommand{\Q}{\mathbb{Q}}
\newcommand{\PP}{\mathbb{P}}
\newcommand{\eps}{\epsilon}
\newcommand{\cN}{\mathcal{N}}
\newcommand{\cM}{\mathcal{M}}
\newcommand{\cA}{\mathcal{A}}
\newcommand{\cB}{\mathcal{B}}
\newcommand{\cD}{\mathcal{D}}
\newcommand{\cH}{\mathcal{H}}
\newcommand{\cL}{\mathcal{L}}
\newcommand{\cF}{\mathcal{F}}
\newcommand{\cP}{\mathfrak{P}}
\newcommand{\cT}{\mathcal{T}}
\newcommand{\cU}{\mathcal{U}}
\newcommand{\cX}{\mathcal{X}}
\newcommand{\iidsim}[0]{\stackrel{\mathrm{iid}}{\sim}}
\newcommand{\NA}[0]{{\tt NA}}
\newcommand{\R}{\mathbb{R}}
\newcommand{\bN}{\mathbb{N}}
\newcommand{\Rp}{\R_+}
\newcommand{\sF}{\mathfrak{F}}

\newcommand{\x}{\mathbf{x}}
\newcommand{\bx}{\mathbf{x}}
\newcommand{\by}{\mathbf{v}}
\newcommand{\xu}{\bar{\mathbf{x}}} %unique design
\newcommand{\yu}{\bar{y}} % averaged observations at \xu
\newcommand{\z}{\mathbf{z}}
\newcommand{\veck}{\mathbf{k}}
\newcommand{\veckN}{\mathbf{k}_N}
\newcommand{\veckn}{\mathbf{k}_n}
\newcommand{\veccn}{\mathbf{c}_n}
\newcommand{\veccN}{\mathbf{c}_N}
\newcommand{\vecX}{\mathbf{X}}
\newcommand{\vecXu}{\bar{\mathbf{X}}} % matrix of unique designs
\newcommand{\vecY}{\mathbf{Y}}
\newcommand{\vecYu}{\bar{\mathbf{Y}}} % Average observations over replicates
\newcommand{\A}{\mathbf{A}}
\newcommand{\An}{\mathbf{A}_n}
\newcommand{\B}{\mathbf{B}}
\newcommand{\C}{\mathbf{C}}
\newcommand{\Cg}{\mathbf{C}_{(g)}}
\newcommand{\CN}{\mathbf{C}_N}
\newcommand{\Cn}{\mathbf{C}_n}
\newcommand{\D}{\mathbf{D}}
\newcommand{\Dn}{\mathbf{D}_n}
\newcommand{\DN}{\mathbf{D}_N}

\newcommand{\Deltan}{\boldsymbol{\Delta}_n}
\newcommand{\I}{\mathbf{I}}
\newcommand{\bT}{\mathbf{T}}
\newcommand{\K}{\mathbf{K}}
\newcommand{\KN}{\mathbf{K}_N}
\newcommand{\Kn}{\mathbf{K}_n}
\newcommand{\SN}{\boldsymbol{\Sigma}_N}
\newcommand{\Sn}{\boldsymbol{\Sigma}_n}

\newcommand{\nug}{\nu_{(g)}}
\newcommand{\hnug}{\hat{\nu}_{(g)}}

\newcommand{\Un}{\boldsymbol{\Upsilon}_n}
\newcommand{\Ug}{\boldsymbol{\Upsilon}_{(g)}}

\newcommand{\U}{\mathbf{U}}
\newcommand{\V}{\mathbf{V}}
\newcommand{\ones}{\bm{1}}
\newcommand{\Lan}{\boldsymbol{\Lambda}_n}
\newcommand{\LaN}{\boldsymbol{\Lambda}_N}
\newcommand{\nun}{\hat{\nu}_n}
\newcommand{\nuN}{\hat{\nu}_N}
\newcommand{\Or}{\mathcal{O}}
\newcommand{\rawTh}{\check{\theta}}
\newcommand{\rawMu}{\check{\mu}}
\newcommand{\rawSigma}{\check{\sigma}}
\newcommand{\modelTh}{\tilde{\theta}}
\newcommand{\modelMu}{\tilde{\mu}}
\newcommand{\modelSigma}{\tilde{\sigma}}
\newcommand{\procMu}{\bar{\mu}}
\newcommand{\procSigma}{\bar{\sigma}}
\newcommand{\procTh}{\bar{\theta}}
\newcommand{\rawU}{\check{u}}
\newcommand{\modelU}{\hat{u}}
\newcommand{\rawPhi}{\check{\varphi}}
\newcommand{\rawRho}{\check{\rho}}
\newcommand{\Var}{\mathbb{V}ar}
\newcommand{\tk}{t_k}
\newcommand{\tkp}{t_{k+1}}
\newcommand{\tK}{t_K}
\newcommand{\tKm}{t_{K-1}}

\newtheorem{lemma}{Lemma}[section]
\newtheorem{Proposition}{Proposition}[section]
\theoremstyle{remark}
\newtheorem{remark}{Remark}[section]

\newcommand{\blu}[1]{\textcolor{blue}{To do: #1}}

%%%%%%%%%%%%%%%%%%%%%%%%%%%%%%%%%%%%%%%%%%%%%%%%%%%%%%%%%%%%%%%%%%%%%%%%%%%%%%%%%%%%%%%%%%%%%%%%%%%%%%%%%%%

\if0\blind
{
\title{\vspace{-1.25cm} A Machine Learning Approach to Adaptive Robust Utility Maximization and Hedging}
\author{Tao Chen\thanks{Department of Statistics and Applied Probability,
University of California Santa Barbara, 1201 Building 479, Santa Barbara, CA 93106-3110}
\and
\and Mike Ludkovski\thanks{Department of Statistics and Applied Probability,
University of California Santa Barbara, 5520 South Hall, Santa Barbara, CA 93106-3110}
}

\date{}

\maketitle
}\fi

\if1\blind
{
  \bigskip
  \bigskip
  \bigskip
  \begin{center}
    {\LARGE\bf A Machine Learning Approach to Adaptive Robust Utility Maximization and Hedging}
\end{center}
  \medskip
  \bigskip
} \fi

\begin{abstract}
We investigate the adaptive robust control framework for portfolio optimization and loss-based hedging under drift and volatility uncertainty. Adaptive robust problems offer many advantages but require handling a double optimization problem (infimum over market measures, supremum over the control) at each instance. Moreover, the underlying Bellman equations are intrinsically multi-dimensional. We propose a novel machine learning approach that solves for the local saddle-point at a chosen set of inputs and then uses a nonparametric (Gaussian process) regression to obtain a functional representation of the value function. Our algorithm resembles control randomization and regression Monte Carlo techniques but also brings multiple innovations, including adaptive experimental design, separate surrogates for optimal control and the local worst-case measure, and computational speed-ups for the sup-inf optimization. Thanks to the new scheme we are able to consider settings that have been previously computationally intractable and provide several new financial insights about learning and optimal trading under unknown market parameters. In particular, we demonstrate the financial advantages of adaptive robust framework compared to adaptive and static robust alternatives.
  \bigskip

\end{abstract}

%\vspace{1cm}

\section{Introduction}
\label{sec:intro}

Stochastic control formulations have been a core tool in financial mathematics for over 40 years. One fundamental challenge in applying stochastic models to practice is the issue of model risk, i.e.~calibrating system dynamics to real life. In this article by ``model risk'' we mean uncertainty about the underlying probability measure
 $\Q$ representing the mis-specification between the assumed dynamics and the true probabilistic structure.  More precisely, the uncertainty concerns some parameters $\theta$, so that we have the parameterization $\Q = \Q^\theta$. Parameter uncertainty arises in any application, since the required model calibration necessarily leaves some residual inference error about the true parameter values. For example, in problems involving risky assets such as stocks, asset return $\mu$ and volatility $\sigma$ are the two main parameters $\theta \equiv (\mu, \sigma)$ driving the investment decisions,  but are notoriously difficult to calibrate.

As one remedy, there have been recently multiple proposals on robust extensions for the underlying stochastic control formulation.  The basic strategy is to take a worst-case view among a collection of potential $\Q^\theta$'s. In the financial context, this conservative perspective stemming from unknown parameters is tempered by the idea of \emph{learning} the dynamics. Historical data about the risky asset evolution can often be used to improve estimates about its drift and/or volatility. Combining these two concepts of robustness and learning, we adopt the \emph{adaptive robust} framework recently proposed in Bielecki et al.~\cite{bielecki2017adaptive}. This paradigm elegantly connects \emph{static robust} approaches that do not allow belief updating with \emph{adaptive} control that fully incorporates the latest estimate within a fixed-parameter setup. Because adaptive control ignores dynamic updating, it  is  time-inconsistent. However, adaptive robust  control also presents formidable numerical challenges as it features an expanded, multi-dimensional state space and non-trivial nonlinear optimization to find the optimal feedback control $\modelU(t,x)$ and the value function $\hat{V}(t,x)$. As a result, closed form solutions are ruled out, and standard PDE methods are not feasible beyond the most simple settings.

In this paper we propose and develop a novel algorithm for adaptive robust control. Our approach belongs to the class of Regression Monte Carlo (RMC) and Control Randomization (CR) strategies \cite{brandt2005simulation,bao2015multi,cong16,han16,zhang2016efficient}. The key idea is to recursively construct a functional approximation $\hat{V}(t,\cdot)$ which is then evaluated over a stochastic (non-gridded) mesh.  However, in the context of robust control, existing techniques  run into the double challenge of: (i) strong path-dependence of optimal $(X_t^*)$ on the control, preventing direct application of Regression Monte Carlo (see Section \ref{sec:hedging}); (ii) the nested optimization due to considering multiple $\Q$'s is computationally intensive and requires proper approximation architecture to extract $u^*(t,x)$. We employ machine learning techniques to overcome the above challenges, by constructing a non-parametric value function approximation \cite{PowellBook,Bertsekas11}. Specifically, we leverage tools from  \emph{statistical emulation} in Stochastic Simulation \cite{ChenRuppert99,BelomestnyKolodko10} (and more broadly the Design and Analysis of Computer Experiments field \cite{Santner2013}) by recasting the task of solving the Bellman equation as a statistical learning problem of fitting a surrogate (i.e.~a statistical model) for $x \mapsto \hat{V}(t,x)$ \cite{risk2016statistical,ludkovski2018simulation}. Our framework resembles recent results for Monte Carlo based solvers of Hamilton-Jacobi-Bellman equations \cite{Pham18deep1,Pham18deep2,guo2019robust}. However, unlike the above references that propose to apply Deep Neural Networks (DNN) for the function approximation step, we rather rely on Gaussian Progress (GP) surrogates. GPs is a core machine learning technique which is well-suited for cases where simulation is expensive and therefore efficient surrogate training is needed (in contrast, DNNs are best suited for high-dimensional problems where maximum flexibility is needed and a lot of training data is available). Related use of GPs is in \cite{balata2019statistical,ludkovski2018simulation}.

Utilizing our numerical algorithm, we provide an extensive investigation into the adaptive robust approach to optimal investment and nonlinear hedging. These are two fundamental problems in financial modeling, yet few fully numerical algorithms have ever been available. Our analysis provides new insights into the interplay between robustifying beliefs, learning and investing/hedging. In particular, we investigate the dependence of the conrols and the terminal wealth on risk- and robustness-level parameters.

The rest of the paper is organized as follows. The rest of Section~\ref{sec:intro} sets up the adaptive robust formulation and lays out the associated discrete-time Bellman recursions. Section~\ref{sec:method} develops our numerical methodology and the Gaussian Process surrogate-based algorithm. The second half of the paper illustrates the algorithm on an Optimal Investment (Section~\ref{sec:utility-max}) and Hedging (Section~\ref{sec:hedging}) case studies. Finally, Section~\ref{sec:enhance} provides discussion on the numerical aspects and concludes with a list of additional enhancements.

\subsection{The Adaptive Robust Bellman Equation}\label{sec:bellman}

% \textbf{Notation:}
Let  $(\Omega, \sF)$ be a measurable space  and $T>0$ a fixed time horizon. Let $\cT=\{\tk: k=0,1,2,\ldots,K,\ t_K=T\}$ be the discrete time set and $\Theta\subset \R^p$   be a non-empty set, which plays the role of the global parameter space, so that for each $\theta \in \Theta$ there is a probability measure $\Q^\theta$ on $(\Omega, \sF)$. We write $\E_\theta$ to denote the expectation operator corresponding to $\Q^\theta$. %By $\PP_{\theta^*}$ we denote the measure generating the true law of $Y$, so that $\theta ^*\in \Theta$ is the (unknown) true parameter. \blu{do we need this?} Model uncertainty implies that $\Theta \supset  \{\theta^*\}$.
We consider a controlled stochastic Markov state process $Y=\{Y^{\vec{u}}_t,\ t\in \cT\}$ taking values in the state space $\mathcal{Y}\subseteq\R^{d}$. We postulate that this process is observed, and we denote by ${\mathbb {F}}=(\cF_t,t\in \cT)$ its natural filtration.
The optimization problem involves the family $\cA$ of admissible feedback control processes, which are ${\mathbb {F}}$--adapted processes $\vec{u}=\{ u_t,\ t\in \cT\}$ defined on $(\Omega,\sF)$ with $u_t$ taking values in a measurable space $\cU$ and being the feedback control applied at time $t$.

Consider the objective of maximizing a reward functional $$C(\vec{u}) := \sum_{k=0}^{K-1} g(\tk, Y^{\vec{u}}_{\tk}, u_{\tk}) + G(Y^{\vec{u}}_{T}),$$
with running rewards $g(t,y,u)$ and terminal reward $G(y)$. In the two financial motivations we have $g(t,y,u) \equiv 0$, but since our algorithm is of independent interest, we continue with this slightly more general setup.
% Above, $Y \equiv (Y^{\vec{u}}_t,\ t\in\cT)$ is the controlled stochastic Markov
% state process taking values in the state space $\mathcal{Y} \subseteq
% \mathbb{R}^d$, $\vec{u} = (u_t)$ with $u_t \in \mathcal{A}$ the feedback control applied at date $t$.
In the classical setup, the dynamics of $Y$ are described through a probability measure $\PP$ and dynamic optimization of ${\cal V}(0,y_0) = \sup_{\vec{u} \in \cA} \E^{\PP}[ C(\vec{u}) | Y_0 = y_0]$ is reduced to the classical dynamic programming Bellman equation
\begin{align}\label{eq:dpe}
{\cal V}(\tk,y) = \sup_{u \in \cU} \left\{ g(\tk,y,u) + \E^\PP_{\tk}\left[ {\cal V}(\tkp, Y^{u}_{\tkp}) \right](y) \right\},
\end{align}
with ${\cal V}(T,y) = G(y)$ and $\E^{\PP}_{\tk}[\cdot](y) \equiv \E^{\PP}[ \cdot | Y_{\tk} = y]$.

%Traditionally, the uncertainty about the $\PP$-dynamics which correspond to some $\theta^*$ has been side-stepped in stochastic control, completely separating the problems of model calibration (inference of $\theta^*$ and subsequent control.

The pervasive challenge of prescribing $\PP$-dynamics implies  model mis-specification for determining $u^*$.
\emph{Robust} control injects consideration of multiple models directly into the optimization step \eqref{eq:dpe} via a generic min-max formulation. Therefore, the single $\PP$ is replaced with many $\Q$'s. In the parametric setup we consider, model uncertainty is indexed via  $\theta \in \Theta$ whose effect is through the dynamics
\begin{align}\label{eq:param-dyn}
  Y^{u}_{\tkp} = T_{(Y)}( \tk, Y_t, u, \theta, \eps_{\tkp}).
\end{align}
The mapping $T_{(Y)}$ captures the joint effect of the initial condition $y$, the control $u$ applied at epoch $\tk$, the parameterized measure $\Q^\theta$, the stochastic shock $\eps_{\tkp}$ and the time-dependence $\tk$. Note that \eqref{eq:param-dyn} postulates that we may re-write the evolution of $Y$ in terms of external, model-independent stochastic shocks; in the common setting we consider, $\eps_{\tk}$ are independent and identically distributed random variables (without loss of generality taken to be standard Gaussian in some $\R^q$).

The robust control setup then includes an intermediate worst-case optimization over a collection $\cM_{k} = \{ \Q^\theta : \theta \in \Theta_{k} \} $, of measures:
\begin{align}\label{eq:robust-dpe}
{\cal V}(\tk,y) &= \sup_{u \in \mathcal{U}} \inf_{\Q \in \cM_{k}} \left\{ g(\tk,y,u) + \E^\Q_{\tk}[ {\cal V}(\tkp, Y^{u}_{\tkp})](y) \right\} \\
 & =\sup_{u \in \mathcal{U}} \inf_{\theta \in \Theta_{k}} \left\{ g(\tk,y,u) + \int_{\R^q} {\cal V}(\tkp, T_{(Y)}(\tk,y,u,\theta, z) )  f_\eps(z) dz \right\}, \label{eq:saddle-point}
\end{align}
where $f_\eps(\cdot)$ is the density of the random variable $\epsilon_{\tkp}$. The expression \eqref{eq:robust-dpe} can be seen as a min-max game between nature that picks $\Q^\theta$ and the controller who counteracts with her control $u=u_{\tk}$. In information-theoretic language, \eqref{eq:robust-dpe} is a game against the adversary who selects a probability measure $\Q^\theta$ among the Knightian uncertainty set $\cM_{k}$ \cite{hansen2006robust}. Note that the adversarial parameter set $\Theta_{k}$ is understood to be dynamic, i.e.~$\Theta_{k} \in \mathcal{F}_{\tk}$ is adapted to the information available by period $\tk$ and in particular can depend on $Y_{\tk}$.

%The key choice then concerns selecting an appropriate uncertainty set $\Theta_t$. A variety of specifications have been proposed, ranging from ``static robustness'' $\theta_t = \theta_0 \forall t$ to dynamic adversarial choice $\Theta_t = \Theta_0$ which corresponds to nature adaptively finding the worst-case parameter setting at each time-step.

The \emph{adaptive robust} control as developed  by Bielecki et al.~\cite{bielecki2017adaptive}, couples the Bayesian updating procedure to the set of measures in \eqref{eq:robust-dpe}.  Specifically, we take
$$
\Theta_{k} = \Theta({\tk}, \procTh_{\tk}),
$$
 re-interpreted as the robustified belief set defined by the latest estimate $\procTh_{\tk}$. As a canonical example, $\Theta({\tk}, \procTh_{\tk})$ is the posterior credible interval at some level $\alpha$ (say the 95\% CI) centered around the unbiased point estimate $\procTh_{\tk}$. This interpolates \emph{dynamic adaptive} control, where $\cM_{k} = \Q^{\procTh_{\tk}} \Leftrightarrow \Theta_{k} = \{ \procTh_{\tk} \} \Leftrightarrow \alpha = 0.5$ (no uncertainty but learning), with the robust versions that do not allow for learning, $\cM_{k} = \{ \Q^\theta : \theta \in \bar{\Theta} \} \Leftrightarrow \Theta_{k} \equiv \bar{\Theta}$ independent of $\tk$. (We should also mention the \emph{myopic adaptive} formulation where one plugs-in $\procTh_{\tk}$ into the static model that treats $\theta$ as a parameter. This is in fact the most common applied usage which myopically separates learning from control.)

 The learning algorithm is meant to reduce the uncertainty about the true probabilistic structure driving $Y$. If one assumes that there exists $\theta^*$ such that $Y_{\tkp} = T_{(Y)}(\tk,y,u,\theta^*, \eps_{\tkp})$ for all $k$ and the learning and estimation are asymptotically perfect, $\procTh_{\tk} \to \theta^*$ as $k \to \infty$ then it is also reasonable to expect $\Theta_{k} \to \{ \theta^*\}$ and the adaptive robust method brings a consistently conservative extension of Bayesian adaptive control. However, the framework can also handle more general situations, such as when $\Theta_{k}$ never converges to a singleton, or when no $\theta^*$ exists. Thus, the adaptive robust can handle non-stationary models where the dynamics change over time (as might be imagined happens in real markets, leading to a process $\theta^*({\tk})$); in such settings $\procTh_{\tk}$ represents the best \emph{current} guess about $\theta^*({\tk})$. So the size of $\Theta_{k}$ need not go to zero.

Returning to the control objective,
since $\procTh_{\tk}$ now affects the inner optimization in \eqref{eq:saddle-point}, it is augmented to the system state
$$
X_t := (Y_t,\procTh_t) \qquad\text{with state space}\quad \cX = \mathcal{Y} \times \Theta \subseteq \R^{d+p}.
$$
The dynamics of $(X_{\tk})=\{X^{\vec{u}}_t,\ t\in\cT\}$ consist of the autonomous dynamics of $(Y_{\tk})$, as well as the filtering equations that govern the evolution of $(\procTh_{\tk})$. We summarize them  by the mapping
$$X_{\tkp} = \bT({\tk},x,u,\theta,\eps_{\tkp})
$$ which aggregates the previous dynamics $T_{(Y)}$ of the state process with the filtering equations for $\procTh_{\tk} \mapsto \procTh_{\tkp}$. Note that those joint dynamics are degenerate since a single noise term $\eps_{\tkp}$ simultaneously drives both $Y_{\tkp}$ and $\procTh_{\tkp}$. As a result, there is a strong dependence between $Y_{\tk}$ and $\procTh_{\tk}$, and more generally $\Theta_{\tk}$, which is to be contrasted with the static robust framework where $Y_{\tk}$ and $\Theta$ are independent. The resulting adaptive robust Bellman equation which is the main object of analysis in this article is then
\begin{align}\label{eq:bellman-recursion}
V({\tk},x) = \sup_{u \in \mathcal{U}} \inf_{ \theta \in \Theta({\tk},\procTh_{\tk})} \E \left[ g({\tk},y,u) + V \left(\tkp, \bT({\tk},x,u,\theta,\eps_{\tkp}) \right) \right].
\end{align}
Note that  $\E$ now simply denotes expectation over $\eps_{\tkp} \in \R^q$ which is the sole stochastic, model-independent component above; the parameter dependence and model risk are encoded in $\bT$. The solution entails applying backward induction over $t=\tKm,\ldots, 0$ to compute $V({\tk},\cdot)$.

\subsection{Solving the Bellman Equation}
To implement the adaptive robust framework, we must solve \eqref{eq:bellman-recursion} which practically reduces to \eqref{eq:saddle-point} with the state-dependent uncertainty set $\Theta({\tk},\cdot)$ as in \eqref{eq:Theta-t}. If we define the operators $\mathscr{F}, \mathscr{M}$ via
\begin{align*}
\mathscr{M}(\tk,x=(y,\procTh))[ \mathscr{F} ] & := \sup_{u \in \cU} \inf_{ \theta \in \Theta({\tk},\procTh)} \mathscr{F}(u,\theta; \tk, x) \\
\mathscr{F}(u, \theta; \tk,x)[V] & := \E \left[ g({\tk},y,u) + V(\tkp, \bT({\tk},x,u,\theta,\eps_{\tkp}) ) \right]
\end{align*}
then the discrete time Bellman recursion is $V({\tk},x) = \mathscr{M}(\tk,x) \circ \mathscr{F}(u, \theta, \tk, x)[V]$. $\mathscr{F}$ is the \emph{propagation operator} that computes the continuation value based on taking conditional expectations of step-ahead value function. $\mathscr{M}$ is the \emph{optimization operator} that computes the optimal control at $(\tk,x)$, recording the optimizers $(\rawTh(\tk, x), \rawU(\tk,x)) = \arg\sup_{u \in \cU} \inf_{ \theta \in \Theta({\tk},\procTh)}$. $\rawTh(\tk,x)$ is the worst-case parameters (i.e.~beliefs) for the model dynamics given $(\tk,y,\procTh)$ and $\rawU(\tk,x)$ is the robustified feedback control in state $x$.
The difference $\rawTh(\tk,(y,\procTh))$ vs $\procTh$ is precisely the precaution taken by the controller against model mis-specification.

Three challenges must be handled to solve \eqref{eq:bellman-recursion}:
\begin{itemize}
  \item Continuous state space: because the state space $\cX$ of $X$ is continuous and  multi-dimensional, some discretization is required to handle it computationally. This discretization however implies that $V(\tk,\cdot)$ (i.e.~$\mathscr{F}(\cdot,\cdot; \tk, x)[V]$) will not be computed for \emph{all} potential $x$, but only for some finite subset. Therefore, space discretization must be accompanied by \emph{interpolation} in order to be able to evaluate  $V(\tk,x)$ for arbitrary $x \in \cX$ in subsequent steps.

  \item Integration: the integral over $\eps_{\tkp}$ in \eqref{eq:saddle-point} requires approximation given that the integrand is not analytically available. This implies that the true conditional operator $\mathbb{E}$ is to be replaced with an approximation $\hat{E}$, i.e.~$\mathscr{F}$ is replaced with $\hat{F}$.

  \item Optimization: finding the optimizers $\rawTh(\tk, x)$ and $\rawU(\tk, x)$ is generally not possible analytically, so another approximation is required to carry out these optimizations numerically. This implies that we replace $\mathscr{M}$ with an approximation $\hat{M}$.
\end{itemize}

In low dimensions, the standard technique for Bellman equations are PDE-based solvers that discretize $\cX$. However, the augmentation of beliefs $\procTh$ to the $Y$-state necessarily leads to multi-dimensional setups, essentially ruling out PDE-based strategies that are practical only in dimension $\le 3$. At the same time, simulation-based strategies from conventional regression Monte Carlo  are also difficult to apply to \eqref{eq:bellman-recursion}:
\begin{itemize}
  \item The control $\vec{u}_t$ intrinsically affects the evolution of the state $(Y^{\vec{u}}_t)$ and hence prevents direct \emph{simulation} of the overall $(X_t)$ as is done in the standard RMC paradigm. %The related Control Randomization method is also hard to apply;

  \item Parametric representation of $\hat{V}(\tk, \cdot)$ (e.g.~in terms of polynomials in $x$) is challenging in dimension $d+p>2$ and brings the concern of overfitting/underfitting;

  \item Because we face the additional step of inner optimization over $\theta$, the representation of $x \mapsto \hat{V}(\tk, x)$ is critical to yield financially reasonable/accurate estimates of $\rawTh(\tk,x)$, which is another limitation of classical parametric approximations;

  \item One usually learns $\hat{V}$ by constructing a grid of $x$-values (akin to finite-difference approaches for PDEs); however such gridding is extremely inefficient for $d+p>2$ and essentially impossible for $d+p>4$.

\end{itemize}

As a result, fully numerical approaches to robust stochastic control have been traditionally viewed as intractable  and there remains a large gap in actual use of robust models for financial modeling. In this article we resolve these challenges, making contributions along two directions. On the algorithmic side, we propose a new, machine-learning-inspired algorithm for \eqref{eq:saddle-point}. The key concept is to employ a non-parametric value function approximation strategy \cite{PowellBook,Bertsekas11}, namely a Gaussian process surrogate. In addition, we borrow concepts from the Design and Analysis of Computer Experiments field \cite{Santner2013}) to construct the underlying stochastic meshes.
This methodology is quite general, and could also be easily modified to tackle other robust control formulations (e.g strong robust), and other contexts (e.g robust optimal stopping). Moreover, it mitigates the issue of scalability to higher dimensions.
 Related extensions would be treated in separate sequels, and we see a lot of further potential for such blending of RMC and machine learning. On the finance side, we present two detailed case studies of using adaptive robust paradigm in optimal investment and loss-based hedging. In particular, the latter problem \eqref{eq:hedge-bellman} was essentially out of reach until now.
 Employing our algorithm, we are able to give a comprehensive investigation of the resulting strategy and its sensitivity to different model parameters. We also demonstrate that the adaptive robust framework outperforms the common (simpler) alternatives of static robust and adaptive control.

\subsection{Motivation: Adaptive Robust Optimal Investment}
The problem of optimal investment in risky financial instruments dates back to Markowitz and Merton and has been extensively studied over the past 50 years. Asset return $\mu$ and volatility $\sigma$ are the two main parameters driving the investment decisions,  but are notoriously difficult to calibrate. Therefore, there is strong interest in methods that explicitly take into account parameter uncertainty, while allowing one to partially learn asset dynamics.

Let $r$ be the constant risk-free interest rate. We assume a fixed time grid $\tk = t_0 + k\Delta t$. The excess log-return of a risky asset $S$ between $\tk$ and $\tkp$ is Gaussian with mean $\mu \Delta t$ and variance $\sigma^2 \Delta t$. Both $\mu$ and $\sigma$  are treated as unknown, yielding the framework of \eqref{eq:robust-dpe} with $\theta \equiv (\mu,\sigma)$. Thus, the dynamics of the controlled wealth process $Y$ under $\Q^{\mu,\sigma}$ are given by
\begin{align}\label{eq:Y-dyn}
  Y_{\tkp}(u) = Y_{\tk}( 1+ r \Delta t + u (e^{\mu \Delta t + \sigma \sqrt{\Delta t} \eps_{\tkp}} - r\Delta t - 1) ),  \qquad \eps_{\tkp} \sim \cN(0,1) \quad\text{i.i.d},
\end{align}
where $u \in [0,1] =: \cU$ is the proportion of portfolio wealth invested in the risky asset. Above we restrict the investor from shorting  $u<0$ or leveraging $u>1$ her risky position.    The objective is to maximize expected utility of terminal wealth
$$\E \left[U(Y^{\vec{u}}_{T}) \right]  \rightarrow \max_{\vec{u}}!$$ for a utility function $U(\cdot)$.
The augmented state $x$ is three-dimensional: $x = (y,\procMu, \procSigma)$ and the Bellman equation for $V(\tk,x)$ is %, where $t$ is now treated as integer (in other words
%we rescale all parameters in the time-scale of the chosen trading frequency $\Delta \tk$, e.g.~$\mu = \mu \Delta t, \sigma = \sigma \sqrt{\Delta t}$) is
\begin{align}\label{eq:inv-1}
  V(t_K,x) &= U(y); \\ \label{eq:inv-2}
  V(\tk,x) &= \sup_{u \in [0,1]} \inf_{ (\mu,\sigma) \in \Theta(\tk,\procMu, \procSigma)} \E \left[ V(\tkp, \bT(\tk,x,u,(\mu,\sigma),\eps_{\tkp}) ) \right]
\end{align}
where the expectation is over the standard Gaussian increment $\eps_{\tkp} \sim \cN(0,1)$.

To specify the transition map $\bT$, we recall the learning of the asset dynamics. Given past observations of $Y$, the respective MLE estimators are denoted  by $\procTh=(\procMu, \procSigma)$ and under $\Q^{\mu,\sigma}$ satisfy for $k=0,1,\ldots$ the recursions~\cite{bielecki2017recursive}
\begin{align}\label{eq:mu-dyn}
  \procMu_{\tkp} & = \frac{ k+1}{k+2} \procMu_{\tk} + \frac{1}{k+2} \big(\mu + \frac{\sigma}{\sqrt{\Delta t}} \eps_{\tkp} \big), \\
  \procSigma^2_{\tkp} & = \frac{ k+1}{k+2} \procSigma_{\tk}^2 + \frac{k+1}{(k+2)^2} \left( (\procMu_{\tk} - \mu)\sqrt{\Delta t} - \sigma \eps_{\tkp} \right)^2. \label{eq:sigma-dyn}
\end{align}
Note that the above are time-dependent, capturing the idea that learning slows over time as the information set grows. As $k$ increases, the weight of the prior $(\procMu_k, \procSigma_k)$ rises and the role of the latest shock $\eps_{\tkp}$ declines.

In sum, the dynamics of  $x_{\tk} = (y_{\tk},\procMu_{\tk}, \procSigma_{\tk})$ under $\Q^{\mu,\sigma}$ are prescribed by \eqref{eq:Y-dyn}-\eqref{eq:mu-dyn}-\eqref{eq:sigma-dyn}, all driven by the 1-D exogenous factor $\eps_{\tkp}$. They are summarized by the map %(recall the substitution $\theta = (\mu,\sigma)$)
\begin{multline}\label{eq:T-optInv}
  \bT( {\tk},(y,\procMu,\procSigma), u, (\mu, \sigma), z) := \Bigl( y(1+r\Delta t+u (e^{\mu \Delta t + \sigma \sqrt{\Delta t} z} - r\Delta t - 1) ), \\ \frac{k+1}{k+2}\procMu + \frac{1}{k+2}( \mu + \frac{\sigma}{\sqrt\Delta t} z), \sqrt{\frac{ k+1}{k+2} \procSigma^2 + \frac{k+1}{(k+2)^2} ( (\procMu - \mu)\sqrt{\Delta t} - \sigma z)^2} \Bigr),
\end{multline}
where $\mu, \sigma$ are treated as external parameters.

Assuming there is a true $\theta^*$,  the uncertainty set $\Theta(\tk, \cdot)$ is described through the (1-$\alpha$)-confidence region for $\theta^*$ for some confidence level $\alpha \in (0,1)$. In the case where only the drift $\mu$ is uncertain and $\sigma$ is given, $\Theta_\alpha(\tk, \procMu_{\tk})$ is an interval centered on $\procMu_{\tk}$,
\begin{align}
 \Theta_\alpha({\tk},\procMu_{\tk}) := \left [\procMu_{\tk} -\frac{\sigma}{\sqrt{ (k+1)\Delta t}}q_{\alpha/2}, \procMu_{\tk} +\frac{\sigma}{\sqrt{ (k+1)\Delta t}}q_{\alpha/2}\right ],
 \end{align}
where  $q_\alpha$ denotes the $\alpha$-quantile  of a standard normal distribution. In the 2-D  case where both $\mu$ and $\sigma$ are uncertain, it is shown in \cite{bielecki2017recursive} that the $(1-\alpha)$-confidence region  for $(\mu^*,\sigma^*)$ at time $\tk$ is an ellipsoid given by
\begin{align}\label{eq:Theta-t}
\Theta_\alpha({\tk},\procMu_{\tk},\procSigma_{\tk}) :=\left\{ (\mu, \sigma) \in\R^2 \ : \ \frac{(k+1)\Delta t}{\procSigma^2_{\tk}}(\mu-\procMu_{\tk})^2 + \frac{k+1}{2 \procSigma^4_{\tk}}(\sigma^2 -\procSigma^2_{\tk})^2\leq \kappa \right\},
\end{align}
where $\kappa$ is the $(1-\alpha)$-quantile of the $\chi^2(2)$ distribution with two degrees of freedom.

When the utility function is of the CRRA power type, $U(y) = \frac{ y^{1-\gamma}}{1-\gamma}$, then the problem \eqref{eq:inv-1}-\eqref{eq:inv-2} is separable with respect to the current wealth $y$: $V({\tk},y,\procMu,\procSigma) = y^{1-\gamma} \tilde{V}({\tk},\procMu, \procSigma^2)$ with the dimension-reduced recursion $\tilde{V}(t_K,\procMu,\procSigma) \equiv \frac{1}{1-\gamma}$, and for $0 \le k \le K-1$
\begin{multline}\label{eq:tilde-V}
  \tilde{V}({\tk},\procMu, \procSigma^2) = \sup_{u \in [0,1]} \inf_{ (\mu,\sigma) \in \Theta({\tk},\procMu, \procSigma)} \E \Big[ (1+r \Delta t + u( e^{\mu \Delta t + \sigma \sqrt{\Delta t} \eps_{\tkp}  } - r\Delta t - 1))^{1-\gamma}  \times  \\ \tilde{V} \big(\tkp, \frac{k+1}{k+2}\procMu + \frac{1}{k+2}( \mu + \frac{\sigma}{ \sqrt{\Delta t}} \eps_{\tkp}), {\frac{ k+1}{k+2} \procSigma^2 + \frac{k+1}{(k+2)^2} ( (\procMu - \mu)\sqrt{\Delta t} - \sigma \eps_{\tkp})^2} \big) \Big].
\end{multline}
We will solve \eqref{eq:tilde-V} in Section \ref{sec:utility-max}. Another case study which also relies on $\theta = (\mu, \sigma)$ is considered in Section \ref{sec:hedging}.

\section{Methodology}\label{sec:method}

Our algorithm is based on value function approximation. This is achieved by evaluating the right hand side of \eqref{eq:bellman-recursion} at a collection of inputs $x^{1:N}_{\tk}$ and then constructing a statistical surrogate for $x \mapsto \hat{V}(\tk, \cdot)$. Thus, we use a machine learning tool, namely Gaussian Process surrogate \cite{rasmu:will:2006}, to learn $\hat{V}(\tk,\cdot)$ based on training data $(x^{1:N}_{\tk}, v^{1:N}_{\tk})$, where $v^n_{\tk} \simeq V(\tk, x^n_{\tk})$, see below. To compute $v^n_{\tk}$ the  expectation in \eqref{eq:tilde-V} is approximated with a quantized sum, and the inner saddle-point is computed locally using a nonlinear optimization algorithm. More precisely we use a numerical quadrature representation $\hat{E}$ of the integral over $\eps_{\tkp}$ and a nested \texttt{fminbnd} call to get $\rawU^n_{\tk}$  which is interpreted as the estimated optimal control at $x^n_{\tk}$. Picking the collection $\cD_k = (x^{1:N}_{\tk})$ is known as experimental design.

%Recall that a surrogate is a functional approximator which can be evaluated via a \texttt{predict} call; for GPs this is done via equation \ref{mean} below.
The overall solver consists of a  \texttt{fit}--\texttt{predict}--\texttt{optimize} loop over time based on the recursion
\begin{align}\label{eq:vn-recursion}
  \hat{V}({\tk},x) = {\sup_{u \in \cU} \inf_{\theta \in \Theta({\tk},\procTh_{\tk})}} \left\{ g({\tk},y,u) + \hat{E} \left[\hat{V}(\tkp, \bT({\tk},x,u,\theta, \cdot) ) \right]  \right\},
\end{align}
where \texttt{predict} is needed to evaluate $\hat{V}(\tkp,\cdot)$ and \texttt{fit} is to construct $\hat{V}(\tk,\cdot)$. For the latter, we evaluate the right hand side of \eqref{eq:vn-recursion} at $x^{1:N}_{\tk}$ which yields the pointwise $v^{1:N})_{\tk}$. Because the latter is based on multiple approximations, including approximating the true $\E[\cdot]$, approximating the true saddle-point, approximating the true $V(\tkp,\cdot)$, we treat them as \emph{noisy} versions of $\hat{V}({\tk},x^n_{\tk})$ and thus the surrogate fitting step includes smoothing, rather than solely interpolation. Similarly, we view $\rawU^n_{\tk} := \arg \sup_{u \in \cU} \inf_\theta \{ g({\tk},y,u) + \hat{E}[ \hat{V}(\tkp, T(x^n_{\tk},u,\theta, \eps) ) ]\}$ as pointwise samples from the feedback control map $\modelU(\tk, x)$ that is also obtained via (an independent) surrogate fitting.

%\begin{enumerate}
%  \item (Assume that the surrogate $\hat{V}(\tkp,\cdot)$ has been fitted)
%
%  \item Select an experimental design $\cD_{k}$ of $N_{k}$ sites $x^n_{\tk}$, $n=1,\ldots, N_{k}$;
%
%  \item Solve the optimization problem \eqref{eq:saddle-point} at each $x^n_{\tk}$, using \texttt{predict} against $\hat{V}({\tkp}, \cdot)$ to evaluate the right-hand-side and a numerical quadrature representation $\hat{E}$ of the integral over $\eps_{\tkp}$. This yields the outputs $v^n_{\tk} := \hat{M}(x^n_{\tk}, \hat{V}({\tkp}, \cdot))$ where $v^n$ is interpreted as the (noisy) estimate of $\hat{V}({\tk}, x^n_{\tk})$ and $\rawU^n_{\tk} := \arg \sup_{u \in \cU} \inf_\theta \{ g({\tk},y,u) + \hat{E}[ \hat{V}(\tkp, T(x^n_{\tk},u,\theta, \eps) ) ]\}$ which is interpreted as the estimated optimal control at $x^n_{\tk}$;
%
%  \item \texttt{Fit} $\hat{V}({\tk},\cdot)$ based on the data $(x^{1:N_{k}}_{\tk}, v^{1:N_{k}}_{\tk})$ and $\modelU({\tk},\cdot)$ based on the data $(x^{1:N_{k}}_{\tk}, \rawU^{1:N_{k}}_{\tk})$;
%
%  \item Goto 1: Start the next recursion for $t_{k-1}$
%\end{enumerate}

%Note that practically the saddle-point is also approximated, and hence the Bellman operator  includes the double approximation of composing $\hat{E}$ and $\hat{M}$ operators. Due to this approximation, we treat the samples $v^n_{\tk} = \hat{M}(x^n_{\tk}, \hat{V}({\tkp}, \cdot))$ as \emph{noisy} signals regarding the true $\hat{V}({\tk}, x^n_{\tk})$, i.e.~the surrogate

%After iterating the above loop for $t=\tKm, \ldots, t_0$, the
Our algorithm returns the fitted surrogates $\{ \hat{V}({\tk}, \cdot) \}$ and $\{ \modelU({\tk}, \cdot) \}$ for each time step. The interpretation of $\hat{V}({\tk},x)$ is as the value function of the underlying adaptive robust control problem. However, for practical purposes we typically wish to know the strategy and resulting utility based on a fixed measure $\Q^{test}$, while \eqref{eq:robust-dpe} assumes that the parameters are dynamically and adversarially drawn from $\Theta_{k}$ at each time-step. To this end, we concentrate on the outputted feedback control surrogate $\modelU({\tk},x)$. To evaluate the resulting utility, we rely on forward Monte Carlo, i.e.~we generate \emph{out-of-sample} forward paths by drawing $Y_{\tkp}$-realizations based on $\Q^{test}$ that induces the controlled trajectories $X^{\modelU}_{0:K}$ and then evaluate the Monte Carlo estimate
$$
\underline{V}^{\Q^*}(0,x_0) = \frac{1}{N'} \sum_{n=1}^{N'} C^n_K.
$$
Observe that $V(0,x_0)$ can also be written as an expectation, but under the path-dependent, ``pasted'' measure $\check{Q}$ which is defined in terms of the worst-case beliefs at each time step
$(\tk, X_{\tk})$ and does not admit a simple interpretation. Indeed, the expectation of $V(\tkp, \bT(\tk, X_{\tk}, \cdot))$ is carried out with respect to $\Q^{\modelTh(\tk,X_{\tk})}$ which depends on $X_{\tk}$. For this reason, the representation $V(0,x_0) = \E^{\check{Q}}[ \sum_k g({\tk},Y^{\modelU}_{\tk},\modelU({\tk}, X_{\tk})) + G(Y^{\modelU}_{t_K})]$ is generally not practically relevant.

\begin{remark}
For an arbitrary $x_*$,
we distinguish between the actual minimizer $\rawU({\tk},x_*)$ that is obtained by  calling \texttt{predict} to evaluate $\hat{V}(\tkp, \bT({\tk},x_*,u,\theta, \eps_{\tkp}) )$ and then optimizing, and the surrogate-based $\modelU({\tk},x_*)$ that is obtained from a statistical prediction. Because optimization is expensive, it is much faster to utilize $\modelU$. We also stress that both $\rawU({\tk},x)$ and $\modelU({\tk},x)$ are defined based on $\hat{V}(\tkp,\cdot)$ and hence will necessarily differ from the true optimal feedback $u^*({\tk},x)$ due to error back-propagation. %This is the price to pay for doing approximate dynamic programming.
\end{remark}

\begin{remark}
  When the optimization over $u$ is over a compact space such as $[0,1]$, an alternative to generic \texttt{fminbnd} or \texttt{optim} gradient-free solvers is to employ direct search over a discrete candidate set such as $\{0, \Delta u, 2\Delta u, \ldots, 1 \}$. This is faster but of course introduces a fixed discretization error relative to the true optimum.  Conversely, the surrogate can output the gradient $\partial \hat{V}/\partial x$ that could be combined with the chain rule to enable gradient-based optimization.
\end{remark}

Before proceeding to discuss the details of operationalizing the above method, we summarize the original proposal of Bielecki et al.~\cite{bielecki2017adaptive} that relied on a Control Randomization (CR) + Regression approach. Bielecki et al.~started with a Monte Carlo-based paradigm that works with trajectories of $(X_{\tk})$ and exploited the fact that the underlying state $(Y_{\tk})$ is an autonomous process. In their approach, one begins by generating $N$ trajectories of $Y_{\tk}$. In addition, one fixes a measure $\Q^0$ and uses the above $N$ trajectories to construct the corresponding filtered posterior estimates $\procTh_{\tk}$ as in \eqref{eq:mu-dyn}-\eqref{eq:sigma-dyn}. Note that $\Q^0$ is picked arbitrarily and serves as a ``baseline'' probability measure like in the CR method \cite{AidLangrene14,kharroubi2014numerical}. Next, the Bellman equation \eqref{eq:robust-dpe} is solved \emph{pathwise} along the above $x^{1:N}_{\tk} = (y_{\tk}^{1:N}, \procTh_{k}^{1:N})$ as follows. Given $({\tk},Y_{\tk},\procMu_{k}, \procSigma_{k})$, one solves for $V({\tk},\cdot)$ by replacing the integral with a weighted sum
\begin{align}\notag
  \hat{E}[\hat{V}(\tkp,x,u,\cdot)] := \int \hat{V}(\tkp, \bT({\tk},x,u,\theta,z) ) f_\eps(z) dz  \\ \approx \sum_{i=1}^{\mathfrak{I}} \hat{V}(\tkp, \bT({\tk},x,u,\theta,\eps^{(i), \mathfrak{I}}) ) w^{(i), \mathfrak{I}},\label{eq:quadrature}
\end{align}
where the Gaussian quadrature recipe specifies the appropriate knots $\eps^{(i), \mathfrak{I}}$ and corresponding quadrature weights $w^{(i), \mathfrak{I}}$. The latter $(\eps^{(i)}, w^{(i)})$ are optimized to minimize a certain global criterion \cite{PagesPrintems05,PagesPhamPrintems04}. Since the noise distribution $\eps_{\tkp} \sim \cN(0,1)$ is fixed throughout, the optimal quadrature recipe is pre-computed and stored offline. In the second sub-step, the one-step-ahead values $\hat{V}(\tkp,\bT(\tk,x,u,\theta,\eps^{(i), \mathfrak{I}}))$ are approximated---since they will not fall on the mesh $\hat{V}(\tkp, x^{1:N}_{\tkp})$---through linear interpolation of the training collection $\{\hat{V}(\tkp,x_{\tkp}^{1:N})\}$.

%In \cite{bielecki2017adaptive} the numerical solution was obtained via discretizing the state space associated with $\widehat \theta$. This was achieved  by generating $I$ sample paths of $Y_t$ under the true dynamics $\theta_*$ and then obtaining the corresponding sample paths of $\widehat{\theta}$ on $[0,T]$. The resulting values $\widehat{\theta}^{1:I}_t$ (viewed as a point cloud in $\R^d$ were used as a stochastic mesh for the space discretization. Then, the recursive equations were numerically solved at each mesh site for the optimal trading strategies $\hat{u}^{1:I} \approx u^*( \widehat \mu_t^{1:I})$.

The above recipe has several limitations that our new proposal overcomes. First, the use of a pre-specified $\Q^0$ leads to a non-adaptive experimental design, which puts $\procTh^{1:N}_{\tk}$ potentially very far from the likely values of $\procTh_{\tk}$ (depending on whether $\Q^0$ is ``good'' or not). By breaking up the path-based construction, we are able to adaptively build better designs, which in particular leverage the structure of $u^*$. Second, linear interpolation of $\hat{V}(\tkp, x^{1:N}_{\tkp})$'s brings three restrictions: (i) it poorly scales in dimension $d$, because the interpolation requires sorting of all the $N$ sites $x^{1:N}_{\tk}$ to identify nearest neighbors; (ii) it implies that $\hat{V}(\tkp,\cdot)$ is piecewise linear which a restrictive approximation architecture (in particular non-smooth which may cause problems when solving for $\rawU$); (iii) it takes the estimated values of $\hat{V}(\tkp,\cdot)$ are taken as \emph{exact}, providing no ability to \emph{smooth} previous approximations. Third, the procedure only outputs the pointwise estimates $\{ \rawU(\tk,x^n_{\tk})\}$; there is no simple way to obtain $\rawU(\tk,x_*)$ for arbitrary $x_*$, essentially ruling out out-of-sample estimation under any other measure than $\Q^0$. Fourth, the approach has no way to avoid \emph{extrapolation} when using \eqref{eq:quadrature}; linear extrapolation is prone to significant errors that are likely to back-propagate. Fifth, pre-simulation under $\Q^0$ requires to use the same mesh size $N$ across all time steps.

\subsection{Gaussian Process Surrogates}
Gaussian process models is a popular choice for flexible statistical models. They allow a consistent treatment of interpolation and regression with very few tunable hyperparameters. The idea of GPs is to view the function to be approximated as a realization of a Gaussian random field with covariance kernel $K(\cdot, \cdot)$. Training the model then reduces to applying the Gaussian conditional equations, i.e.~evaluating the distribution of the random field given the data. The Gaussian structure implies that the conditional process is still Gaussian, so its distribution is summarized through the posterior mean $m_*(\cdot)$ and posterior (co-)variance $s_*(\cdot, \cdot)$, interpreted as the model-predicted approximation $\hat{V}(\tk,\cdot)$ and the corresponding posterior uncertainty or ``standard error''.

Intuitively, the GP offers a methodology to put a smooth surface through the data labelled generically as $(\bx, \by) \equiv (x^{1:N}, v^{1:N})$. For numerical stability, we introduce the ``nugget'' $\eta^2$, so that rather than exactly interpolating the given $v^n$'s, the surrogate also smoothes. Namely, the nugget is equivalent to assuming
 $$v^n_{\tk} = V(\tk, x^n_{\tk}) + \eps^n \qquad \text{with}\quad \eps^n \sim \mathcal{N}(0,\eta^2).$$
 In our context, adding $\eps^n$'s is justified via the errors coming from quantization and numerical optimization and we use $\eta = 10^{-5}$. The predicted $\hat{V}(\tk, x_*) \equiv m_*(x_*)$ at  any input $x_*$ is given by
\begin{align}
{m}_*(x_*) &=  \mathbf{k}(x_*)[\mathbf{K}+\eta^2\mathbf{I}]^{-1}\by,  \label{mean}\\
\intertext{with corresponding posterior covariance $s_*(x_*, x_*')$}
s_*(x_*,x_*') &=  K(x_*,x_*')-\mathbf{k}(x_*) [\mathbf{K}+\eta^2 \mathbf{I}]^{-1} \mathbf{k}(x_*'), \label{cov}
\end{align}
with the $N \times 1$ vector $\mathbf{k}(x_*)$ and $N \times N$ matrix $\mathbf{K}$ defined by
\begin{align}
\mathbf{k}(x_*) := K(x_*,\bx) &= [K(x_*,x^1), \ldots ,K(x_*,x^N)],\quad \text{and}\quad \mathbf{K}_{i,j} := K(x^i,x^j).
\end{align}
Above $\mathbf{I}$ is a $N \times N$ identity matrix and $\eta^2 \mathbf{I}$ is the noise matrix.
%and the diagonal of $K_\vartheta$ is $K_\vartheta(x,x) = \tau^2 + \sigma^2$.
%
%For our purposes, we may concentrate on the posterior mean $m_*(\cdot)$ which is viewed as the surrogate $\hat{V}(t,x) \equiv m_*(x)$.

%\subsection{GP Fitting and Implementation Details}
To construct a GP surrogate requires selecting a kernel \emph{family} and the corresponding hyper-parameters. The latter step is typically done through Maximum Likelihood Estimation, i.e.~a nonlinear optimization problem involving the respective likelihood of observations.  A common kernel family we use is the Matern-5/2,
\begin{equation}
\label{eq:matern52}
K_\vartheta(x,x'):= \tau^{2} \prod_{i=1}^{d+p} \left(1+ \frac{\sqrt{5}r}{\rho_i} + \frac{5r^{2}}{3\rho_i^{2}}\right)\exp\left(-\frac{\sqrt{5}r}{\rho_i}\right), \quad r:=|x_{i} - x'_{i}|,
\end{equation}
where the lengthscales $\rho_i$ and the process variance $\tau^2$ are estimated via maximum likelihood. Specifically, we use the GPML \texttt{Matlab} package to obtain the MLE for the hyperparameters $\vartheta := (\eta, \tau, \rho_{1:(d+p)})$. Another kernel is the squared-exponential
\begin{equation}
\label{eq:sq-exp}
K_\vartheta(x,x'):= \tau^{2} \prod_i \exp\left(-\frac{(x_{i} - x'_{i})^2}{2 \rho^2_i}\right).
\end{equation}
%which yields $C^\infty$ smooth $m_*(\cdot)$.

The GP predictive surface $x_* \mapsto m_*(x_*)$ is akin to kernel regression in the sense that the prediction \eqref{mean} at $x_*$ is always a weighted average of the $v^n$'s, with the weights driven by the spatial covariance structure encoded in $K(\cdot, \cdot)$. In particular, if the spatial correlation decays quickly, the surface will tend to be more ``bumpy'', while for very strong spatial correlation the surface will be nearly flat/linear. The differentiability of the predictive surface $m_*(\cdot)$ is driven by the properties of $K(\cdot,\cdot)$. Under the Matern-5/2 choice \eqref{eq:matern52}, the resulting $m_* \in \mathcal{C}^2$ is twice-differentiable, while for \eqref{eq:sq-exp} it is ${\cal C}^\infty$.  In both families, the lengthscales $\rho_i$ determine the spatial ``wiggliness'' of the fitted surface in the respective coordinate. Note that different lengthscales (anisotropy) in different coordinates allow for $\hat{V}$ to be, say, more flat in $S$, but more flexible in $\procMu$.

%\subsection{Spline Surrogates}

For the propagation operator $\mathscr{F}$  showing in the Bellman equation we proceed as in \eqref{eq:quadrature},
approximating the integral over a standard normal distribution by a $\mathfrak{I}$-points optimal quantizer $\hat{E}$, i.e.~a discrete sum over $\mathfrak{I}$ quadrature points $\eps^{(i),\mathfrak{I}}, i=1,\ldots,\mathfrak{I}$.

\begin{remark}
  In CR-RMC \cite{AidLangrene14,kharroubi2014numerical}, solving the Bellman equation
  $$
  v(\tk,x) = \sup_{u \in \cU} \E \left[ v(\tkp,X^{u,x}_{\tkp}) \right]
  $$
  is handled by defining the q-value $q(\tk,x,u) := \E[ v(\tkp,X^{x,u}_{\tkp})]$ across $u \in \cU$, building a surrogate $(x,u) \mapsto \hat{q}(\tk,x,u)$ and then
  finding $\rawU(\tk,x)$ as the analytic maximizer of $u \mapsto \hat{q}(\tk,x,u)$. In contrast, in our approach, we first fix $x$ and maximize $\mathscr{F}$ at $x$ to obtain $\rawU(\tk,x)$, then fit a surrogate $\modelU(\tk,x)$, effectively reversing the order of optimization and emulation. %This is similar to the RMC distinction between regress-later (project, then take conditional expectation) and regress-now (take conditional expectation, then project), cf.~\cite{balata17,balata18}.
\end{remark}

\subsection{Experimental Design}
%\textbf{Simulation Design:}
The quality of the regression estimates is critically linked to the choice of the simulation design $\cD_{k} = x^{1:N}_{\tk}$. Specifically, accuracy of the GP-based $\hat{V}(\tk,x_*)$ at some given input $x_*$ is directly related to the \emph{density} of the design $\cD_{k}$ around $x_*$. This is the localization property of the GP and matches the intuition of constructing an interpolant for the training data $(x^{1:N}, v^{1:N})$: the interpolated prediction would be good in the neighborhoods of $x^n$'s, and progressively worse far away. In particular, as $x_*$ gets far from the training inputs $\bx$, the GP prediction $m_*(x_*)$ is driven by the chosen asymptotes of $\hat{V}(\tk,\cdot)$ reflecting extreme extrapolation. Specifically, for GPs we have that $m_*(x_*) \to m(x_*)$ reverts to its prior mean.

Based on the above discussion, the design $\cD_{k}$ should concentrate on the \emph{region of interest}. In our context, the latter is driven by the pairs $(\procMu_{\tk}, \procSigma_{\tk})$ that are likely to be encountered by the controller at step $k$. In addition, since the key output we seek is the investment strategy $u^*({\tk},x) \in [0,1]$, the hardest learning task is to identify the optimal control when the constraints are not binding. Consequently, the region of interest is $R = \{ x: u^*({\tk},x) \in (0,1), p({\tk},x | x_0) > \underline{p} \}$.

\begin{figure}[ht]
  \centering
    \includegraphics[width=0.45\textwidth]{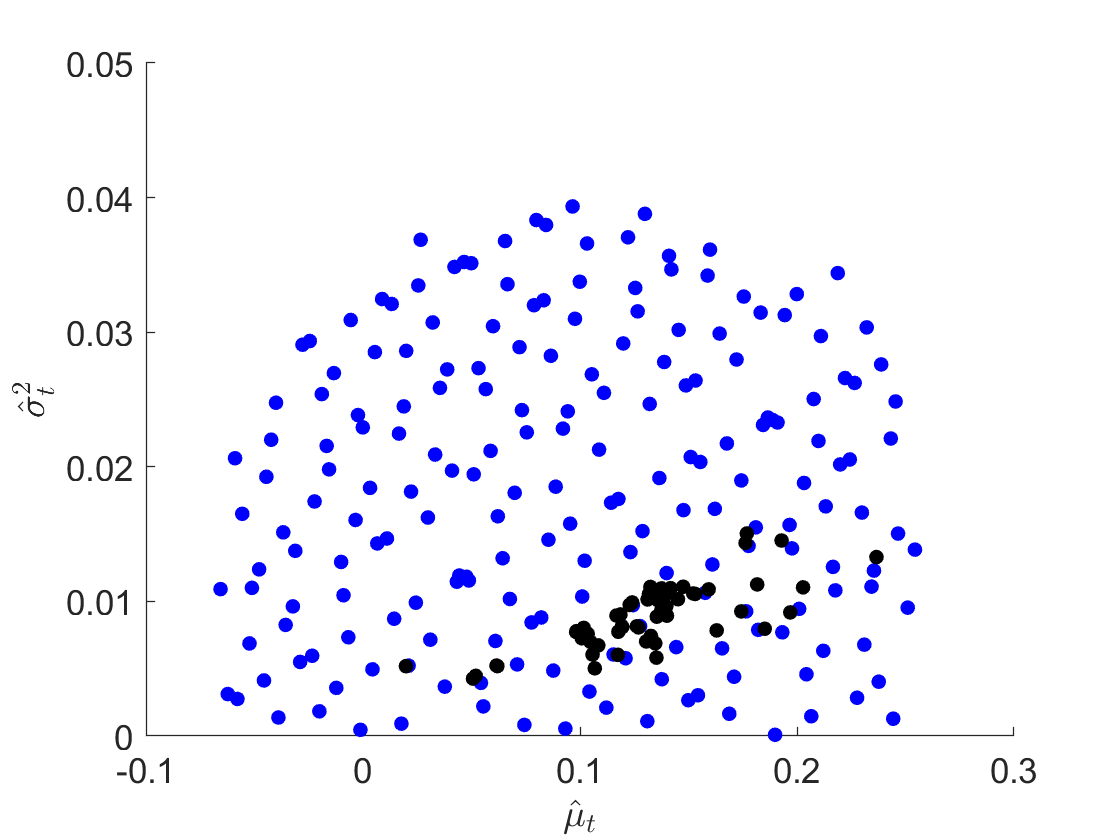}  \includegraphics[width=0.45\textwidth]{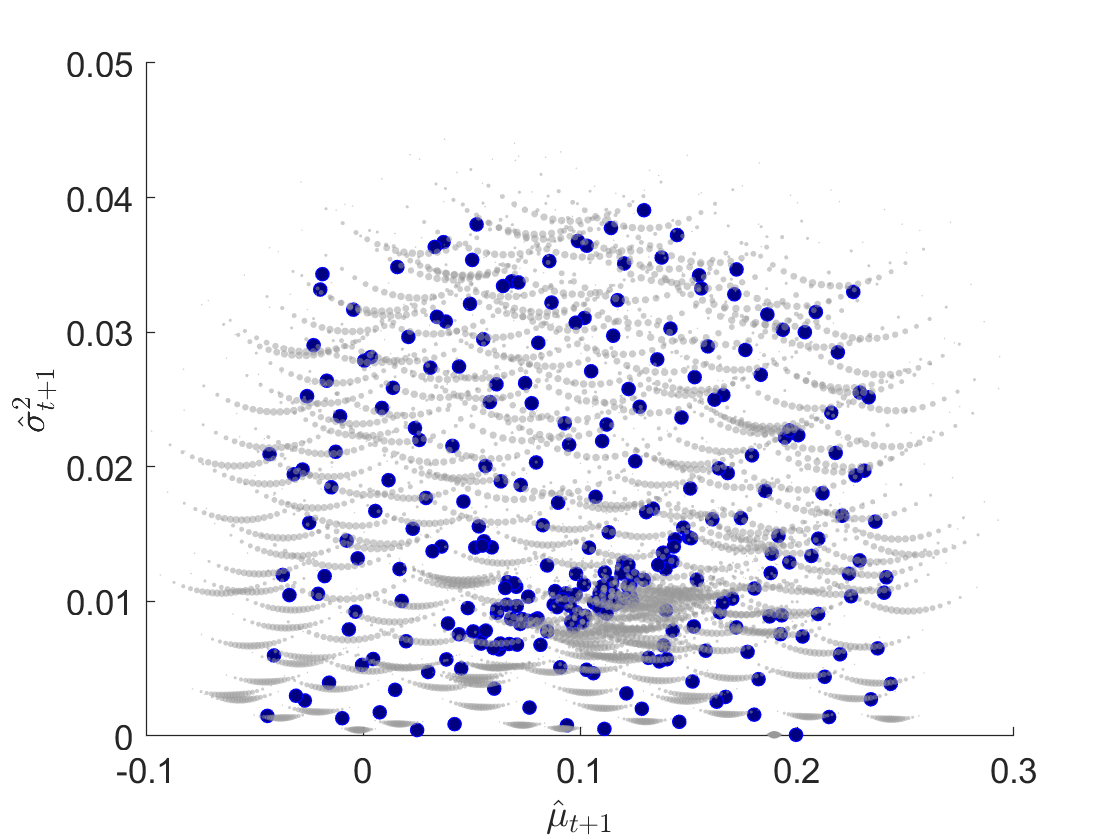}
 \caption{Constructed simulation design $\cD = \{ \procMu^{1:N}, \procSigma^{1:N} \}$ for the portfolio optimization case study. There are a total of 256 sites: Sobol QMC based (200 blue); adaptive based on $\modelU_{\tkp}$ along pilot paths (56 black). The right panel shows the predictive locations $V(\tkp,\cdot)$ used from $\cD_{k}$, with the blue dots indicated the respective $\cD_{k+1}$ to highlight inter- and extra-polation.  \label{fig:sim-design}}
\end{figure}

Another problem that requires attention is extrapolation in estimating the surrogate $\hat{V}({\tk},\cdot)$. Extrapolation tends to lead to large estimation errors in approximating the true conditional expectation operator $\mathbb{E}$.  With this in mind, we need to select the design $\cD_{k}$ so that extrapolation is minimized. In order to compute $\hat{V}({\tk},\procTh_{\tk}^n)$,  we must consider $\hat{E}[\hat{V}(\tkp,\mathbf{T}({\tk},x_{\tk}^n,u,\theta,\epsilon_{\tkp}))]$ at any $\theta$ within $\Theta({\tk},\procTh_{\tk}^n)$. This, in turn, means that we need good estimation of $\hat{V}(\tkp,\mathbf{T}({\tk},x_{\tk}^n,u,\theta,\epsilon^{(i)}))$, $i=1,\ldots,\mathfrak{I}$, where recall that $\epsilon^{(i)}$ comes from Gaussian quantization. Figure \ref{fig:sim-design} displays the relationship between $\cD_{k}$ and the set of all next-step locations $\{ x_* \in \mathbf{T}({\tk},x_{\tk}^n,u,\rawTh(\tk,\procTh^n_k),\epsilon^{(i)})), n=1,\ldots, N; \epsilon^{(i)} \in Quantizer_{\mathfrak{I}} \}$. These are the locations where we must evaluate $\hat{V}(\tkp, x)$ and hence highlights to what extent the method requires extrapolation. On the right panel we compare all needed predictive sites $\theta_*$ against the $\procTh^{1:N}_{\tkp}$ used for training $\hat{V}(\tkp,\cdot)$.

Given a region of interest $\hat{R}_{\tk}$, there are two main concepts to generate a training set $\cD_{k}$ that targets $\hat{R}_{\tk}$. The first idea is  a ``density-based'' approach that aims to make $\cD_{k}$ mimic the distribution of $X_{\tk}$ under $\Q^{test}$. Recall that $X_{\tk}$ includes the exogenous $Y_{\tk}$ and the path-dependent $\procTh_{\tk}$. One method to handle this path-dependency is based on Control Randomization~\cite{AidLangrene14,kharroubi2014numerical} which was also the motivation for generating the path-simulations in Bielecki et al.~\cite{bielecki2017adaptive}. We first pick some $\Q^0$ and generate paths under $\Q^0$. In other words, we pick some ``true'' model parameters  (e.g.~using the prior mean estimates $\procTh_{t_0}$) that are used to generate $Y_{\tk},\procTh_{\tk}$. A collection of $N$ such independent ``pilot'' paths can be then used to set $\cD_{k} := \{ (y^n_{\tk}, \procTh^n_{\tk}) \}$. %As long as $\Q^0$ is close to  $\Q^*$ this will properly explore the state space of $X^*_{\tk}$.

The second concept is to apply space-filling schemes which aim to yield a ``uniform'' sample from the region of interest. This approach is similar to classical maximin/max-entropy experimental designs in statistics. Space filling is achieved by specifying an input domain (usually some polyhedral or rectangular bounding box) and then a space-filling sequence. A popular choice are Quasi Monte Carlo (QMC) sequences that are used as a variance reduction tool to ensure that $\cD_{k}$ does not contain any ``holes'' that might degrade estimation in their neighborhoods. A QMC sequence offers a deterministic way to sequentially fill the unit hypercube and can be straightforwardly scaled and clipped to fill any polyhedral domain. For example, Sobol QMC sequences offer a discrete $\cD_{k}$ of any size $N$. %mitigates extrapolation errors or lack of knowledge about $\Q^*$ and

% In particular, two concerns arise: if the design concentrates mostly around some region $R$, then the estimates would be accurate around $R$ and less accurate far from it. Second, if $x$ is far from $x_k$, then prediction boils down to extrapolation, i.e.~rather than being data-driven, the answer is effectively driven by the chosen asymptotes of $\hat{V}$ as $x \to \infty$. For example, for polynomial regressions, $\hat{V} \to \pm \infty$ outside the training set.

% Idea 1: guess $R_k = R_{k+1}$. Idea 2: space-filling. Idea 3: guess $Q^0$ and generate forward paths.

In Figure \ref{fig:sim-design}, we show an adaptive ``mixture'' design $\cD_{k}$ that we utilize for the portfolio optimization problem with unknown $(\mu,\sigma)$. It is obtained in three steps. In the first step, we simulate $N''=250$ forward paths under a pre-specified $\Q^0$ that yield $(\procMu''_{\tk}, \procSigma''_{\tk})$. We then utilize these ``pilot'' paths to obtain a convex hull $R'_k$. In the second step, we utilize the Sobol QMC sequence to fill the respective $R'_{\tk}$ with $N'=100$ sites. Finally, in the third step, we augment with another $\tilde{N}=50$ sites based on $R_{\tk}= \{ x_{\tkp}^n : \modelU(\tkp, x^n_{\tkp}) \in (0,1) \}$, i.e.~add new sites where the control was non-trivial in the previous time-step.

%In addition,  we utilize \emph{sensible} fall-back rules for respecting the constraint $u^* \in [0,1]$ during extrapolation. %For the region of interest $R_k$, we utilize later time-steps to construct the design at $t_k$, namely $R_k = \{ x_{k+1}^n : \modelU(k+1, x^n_{k+1}) \in (0,1) \}$, i.e.~pick the new input space as the region of interest from the previous step $k+1$.

%implemented a simple variant of the second and third ideas by letting $\cD_k$ be a ``mixture'' design. Specifically, we selected 100 sites based on a Sobol QMC sequence which is space-filling in a pre-specified rectangular domain. We then select another 150 sites by sampling $(\procMu_t, \procSigma_t)$ from a pre-specified law $\Q^0$ (for simplicity we have used $\Q^0 = Q^*$ to match the true model used for testing, which clearly would not be available in real-life context).

%A more advanced method would be to have an adaptive design, for example by growing $\cD_k$ on the fly (sequential design) as we discover what states are most relevant.

\emph{Remark:} Given the complexity of \eqref{eq:bellman-recursion} which requires solving a min-max problem at each design site $x^n_{\tk}$, it is not computationally feasible to have thousands of $x^n$'s. Already for a rather small design like in Figure~\ref{fig:sim-design}, we are looking at more than a million of \texttt{predict} calls to the GP surrogate \emph{each time step} ($\mathfrak{I}=100$ for each integral, multiplied by about 60 optimization steps to converge to the min-max optimal value, multiplied by $|\cD|=N= 256$ design sites). This is a principal reason why we advocate Gaussian Process surrogates which excel in learning from small-to-moderate designs. Note that GPs are still computationally intensive and their speed is very sensitive to $N$, which is the reason for all the above care in building a good experimental design. While a linear surrogate that utilizes ordinary least squares regression with a pre-specified set of basis functions is much faster, its rigid structure performs terribly for learning $u^*$; in our experience success of the algorithm hinges on having enough degrees of freedom. Other alternatives for value function approximation besides GPs are left for future research.

\subsection{Algorithm}\label{sec:algo}

Algorithm \ref{algo} summarizes our procedure. It further makes a slight extension to allow for time-dependent design size $N_k$. The second Algorithm \ref{algo:fwd} evaluates the performance and strategy of the controller under some test measure $\Q^{test}$ utilizing a forward Monte Carlo average over $N'$ fresh paths. The actual simulations are under $\Q^{\theta^{*,n}}$ where $\theta^{*,n}$ can vary across paths or even across time. For example, if we take $\theta^{*,n} \sim \Theta_0$ uniformly, then we mimic the static robust setup. The resulting out-of-sample estimator  $\underline{V}$ is guaranteed to be a lower bound for $V(0,x_0)$ since it is directly based on the suboptimal strategy $\modelU(\cdot, \cdot)$ and its only error is the Monte Carlo averaging due to the finiteness of $N'$.

\begin{algorithm}[ht]
\caption{Backward Recursion to learn $\hat{V}$ and $\modelU$. \label{algo}}
\begin{algorithmic}[1]
\REQUIRE Design sizes $N_k$
\STATE Initialize with terminal condition $\hat{V}(K,x,\theta) = G(y)$ $\forall x$.
\FOR{ $k = K-1, \ldots, 1$}
\STATE Create a design $\cD_{k} = (Y^n_{\tk}, \procTh^n_{\tk}) \equiv x_{\tk}^{1:N_{k}}$ that will be used to estimate $\hat{V}(t_{\tk},\cdot)$.
\FOR{$n=1,2,\ldots, N_{k}$}
\STATE Using an optimal Gaussian quantizer with $\mathfrak{I}$ terms set
  \begin{align}\label{eq:f2}
  \mathscr{F}_2(u, x^n_{\tk}) := \inf_{ \theta \in \Theta({\tk},\procTh^n_{\tk})} \sum_{i=1}^{\mathfrak{I}} \hat{V}(\tkp, \bT({\tk},x^n_{\tk},u,\theta,\eps^{(i), \mathfrak{I}}) ) w^{(i), \mathfrak{I}}.
  \end{align}
\STATE Let \begin{align}\label{eq:yn}
  v^n_{\tk} := \sup_{u \in {\cU}} \mathscr{F}_2(u; x^n_{\tk}).
   \end{align}
\STATE Record the estimated optimal control $\rawU^n_{\tk} \leftarrow \arg\sup \mathscr{F}_2(u; x^n_{\tk})$.
\STATE Record the driving worst-case parameters corresponding to $\rawU^n_{\tk}$: $$\rawTh^n_{\tk} \leftarrow  \arg\inf_{\theta} \sum_{i=1}^{\mathfrak{I}} \hat{V}(\tkp, \bT({\tk},x^n_{\tk},\rawU^n_{\tk},\theta,\eps^{(i), \mathfrak{I}}) ) w^{(i), \mathfrak{I}}.$$
\ENDFOR

\STATE Build a GP model $\hat{V}({\tk},\cdot)$ for the link between $(x^{1:N_k}_{\tk})$ and $(v^{1:N_k}_{\tk})$; this is the functional representation of the value function at step $k$.
\STATE Build a GP model $\modelU({\tk},\cdot)$ for the link between $(x^{1:N_k}_{\tk})$ and $(\rawU^{1:N_k}_{\tk})$; this is a (separate) functional representation of the optimal adaptive robust feedback control map %Build a GP model $\hat{\Theta}(t,\cdot)$ for the link between $(x^n_{k})$ and $(\check{x}^n_k)$.
\ENDFOR
\end{algorithmic}
\end{algorithm}

\begin{algorithm}[ht]
\caption{Forward Monte Carlo to evaluate performance of  $u^*$ under $\Q^{test}$ \label{algo:fwd}}
\begin{algorithmic}[1]
\REQUIRE No.~of simulations $N'$, no.~of time-steps $K$, initial $(X_{t_0}, \procTh_{t_0})$.
\STATE Initialize $\procTh^n_{t_0}  \leftarrow \procTh_{t_0}, x^n_{t_0} \leftarrow x_{t_0}, \  n=1,\ldots N'$
\STATE Set $\theta^{*,n}$ which is the parameter set (possibly time-dependent) for the $n$-th forward path (so simulated under $\Q^{\theta^{*,n}}$).
\FOR{$k=0,\ldots,K-1$}
\STATE Draw i.i.d. $\eps^n_{\tkp} \sim {\cal N}(0,1)$, $n=1,\ldots, N'$
\STATE Using the GP surrogate compute the control $u^n_{\tk} \leftarrow \modelU(t_{\tk}, x^n_{\tk})$.
\STATE (Optional) Evaluate (using prediction from the surrogate or a direct evaluation of the optimizer) the worst-case parameters $\modelTh^n_{\tk}$ which are functions of $x^n_{\tk}$;
\STATE Compute the realized payoff $g(t_k, y^n_{\tk}, u^n_{\tk} )$; update the cumulative $C^n_{k} \leftarrow \sum_{\ell=1}^k g(t_\ell,y^n_{t_\ell},u^n_{t_\ell})$;
\STATE Update the states according to $x^n_{\tkp} \leftarrow \bT \left({\tkp},x^n_{\tk},u^n_{\tk}, \theta^{*,n}, \eps^n_{\tkp} \right).$
\ENDFOR
\STATE Return $\underline{V}(0,x_0) := \frac{1}{N'} \sum_{n=1}^{N'} C^n_{K}$
\end{algorithmic}
\end{algorithm}

%We initialize with a given starting condition $\procMu^n_{t_0}  = \procMu_{t_0}, \procSigma^n_{t_0} = \procSigma_{t_0}, x^n_{t_0} = x_{t_0} \ \forall n=1,\ldots N'$. Now repeat over $k=0,\ldots K-1$:
%
%\begin{enumerate}
%  \item Draw i.i.d pairs $Z_{\ell}|Z_{\ell-1} \sim \Q^*$ which are identified via $\eps^n_{\tk}$, $n=1,\ldots, N'$;
%
%  \item Using the GP surrogate compute and apply the control $u^n_{\tk} = \modelU(t_{\tk}, x^n_{\tk})$. Further evaluate (using prediction from the surrogate or a direct evaluation of the optimizer) the worst-case parameters $\modelMu^n_{\tk}, \modelSigma^n_{\tk}$ which are functions of $x^n_{\tk}$;
%
%  \item Compute the realized payoff $g(t_k, x^n_{\tk}, \modelU(t_k, x^n_{\tk}) )$ and update the cumulative $C^n_{k} = \sum_{\ell=1}^k g(t_\ell,y^n_{t_\ell},u^n_{t_\ell})$;
%
%  \item Update the states according to $x^n_{\tk} = \bT \left({\tk},x^n_{t_{k-1}},u^n_{\tk},\eps^n_{\tk}; \procMu^n_{t_{k-1}}, \procSigma^n_{t_{k-1}} \right).$
%\end{enumerate}
%The final answer is the average $\underline{V}(0,x_0) := \frac{1}{N'} \sum_{n=1}^{N'} C^n_{K}$. %Note that a different (statistically worse) estimate is simply $\hat{V}(0,x_0)$. In contrast to the latter which is likely to be biased,
% (which disappears as $N' \to \infty$. One common trick is to use $\hat{V}(0,x_0)$ as an empirical upper bound and $\underline{V}(0,x_0)$ as the empirical lower bound for the final answer.

To summarize, the proposed framework organically integrates the construction of $\hat{V}$ and $\modelU$ with the paradigm of Design and Analysis of Computer Experiments that views the intermediate step as inference of an \emph{expensive black-box} function. The adaptive $\cD$'s target the learning of the optimal control, while the GP setup captures the spatial borrowing of information  to predict $\modelU(\tk,x)$ without directly optimizing.
Their combination improves efficiency, scalability, and interpretability, and should be contrasted to the conventional implementation where $\cD$ is a grid and the approximation architecture $\cH$ is a linear interpolant (piecewise linear $\hat{V}$). Our algorithm is also highly modular, allowing additional ways of approaching the aforementioned sub-problems of numerical integration and numerical optimization. Intuitively, the computation revolves around repeated optimization, and our framework achieves substantial gains by leveraging already obtained solutions of similar optimization problems, in analogy to \emph{parametric optimization}. To our knowledge, we are the first article to propose this strategy for robust stochastic control.

\begin{remark}
In the case studies below we have $\cU = [0,1]$ and the GP surrogate $\modelU$ tends to over-smooth in regions where $\rawU({\tk},x) \simeq 0$ or $\rawU({\tk},x) \simeq 1$ as it does not like flat responses. We regularize the constraint $u \in [0,1]$  by carrying out optimization of $\mathscr{F}_2(u,x)$ over a larger domain $u \in \tilde{\cU} \supseteq \cU$. The relaxed maximizer $\rawU(x)$ is then fed into the GP $\modelU$ representation of the optimal control, and the resulting prediction is \emph{projected} as $\mathfrak{P}_\cU (\modelU) \in [0,1]$. This regularization makes sure that the data $\rawU^n_{\tk}$ used to build the surrogate are smooth and do not have the non-smooth cut-offs at 0 and 1. In the first example below we take $\tilde{\cU} = [-0.2, 1.2]$.

\end{remark}

\subsection{Algorithm Stability}
The overall algorithm complexity can be broken down into the overhead of fitting the surrogates and the cost of calling \texttt{predict} during the Bellman recursions. With a GP surrogate, the complexity of model fitting is $\Or(N_k^3)$ at each time step, so $\Or( \sum_k N^3_k)$ total. Each GP prediction requires $\Or( N_k^2)$ effort and is employed $\mathfrak{I}$ times during each evaluation of the conditional expectation $\hat{E}$ and in turn is called $N_{optim}$ times (which is state-dependent) by the nested optimizer that is solving for $\rawU$. Overall, we obtain $\Or( \mathfrak{I} \times N_{optim} \times \sum_k (N_k^3))$ complexity for Algorithm \ref{algo}. Similarly we have $\Or( \sum_k N_k^2 \times N')$ complexity for generating $N'$ forward paths to obtain the out-of-sample evaluation of $\underline{V}$ in Algorithm \ref{algo:fwd}.

An important aspect is that the proposed algorithm is non-deterministic since the outputted $\hat{V}$ and $\modelU$ depend on the underlying random samples of $\eps^{1:N_k}_{\tk}$. Recall that $\eps^n_{\tk}$ affects the $v^n_k$ that is used for the surrogate fitting. In the GP context, this will modify both the GP hyperparameters, as well as the actual prediction since $\by$ shows up on the RHS of \eqref{mean}. Therefore, running the algorithm twice will generate slightly different fits. The respective variance would be amplified if the simulation designs $\cD_{k}$ are also random (e.g.~dependent on the pilot paths as in our setup) since then also the $x^n$'s would vary from run to algorithm run. Yet another increase in algorithm variance would occur if also the conditional expectation is approximated via Monte Carlo rather than Gaussian quadrature.

\begin{figure}[ht]
  \centering
    \includegraphics[width=0.45\textwidth]{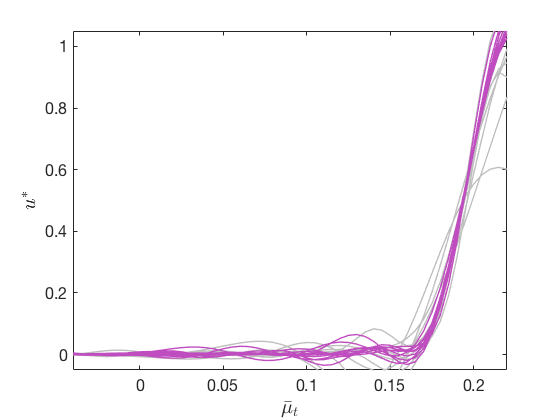}
    \includegraphics[width=0.45\textwidth]{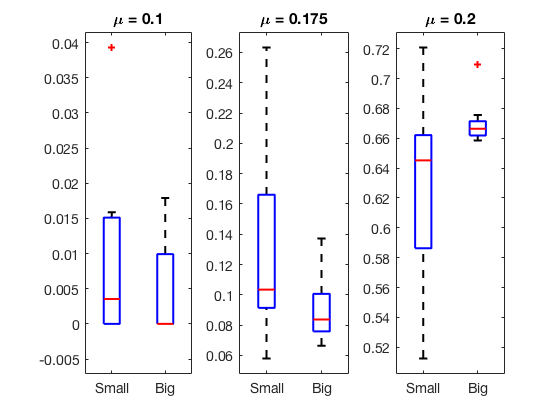}%{predictedControl_MCerror.png}  %\includegraphics[width=0.45\textwidth]{RegressionPictures/staticAComp-v2.png}
 \caption{Comparison of Monte Carlo standard errors as a function of design size $N \equiv N_k \forall k$. Left: Ten empirical predictions $\modelU(\tk,\procMu_{\tk}, \procSigma_{\tk})$ at ${\tk}=0.8$, $\procSigma^2_{\tk} = 0.01$ as a function of $\bar{\mu}_{\tk}$ to illustrate the intrinsic fluctuations of our algorithm. We show both low-budget (in purple) and high-budget (in grey) runs. Right: boxplot for $\modelU(\tk, \procMu_{\tk},0.01)$ at three different $\procMu_{\tk}$. \label{fig:MC-error}}
\end{figure}

In this vein, the design size $N_k$ plays a double role of (i) decreasing the above macro-replication variance through Law of Large Numbers with respect to $\eps^n_{\tk}$; (ii) increasing accuracy (i.e.~reducing bias) by raising the fidelity of $\modelU$ and $\hat{V}$ (i.e.~reducing functional approximation error between the surrogate and the true minimizers/value function). Figure \ref{fig:MC-error} displays a boxplot of $\modelU$ across algorithm macro-replications when using $N_k \equiv 100$ vs $N_k \equiv 250$. As expected, a larger design reduces overall variance, as well as the bias. To illustrate the role of another tuning parameter, we show the impact of the quantization level $\mathfrak{I}$ on the accuracy of the results.

%\blu{Run algorithm with two different budgets to see impact on accuracy. Show GP posterior variance for 2 simulation designs of different size $|\cD| = 100$ vs $|\cD| = 250$. Figure \ref{fig:MC-error} shows a boxplot of $u^*$ at reference state to  show algorithm stability across macro-replications.}

\section{Adaptive Robust Utility Maximization}\label{sec:utility-max}

In this section we illustrate our algorithm from Section \ref{sec:algo} for the task of portfolio optimization under uncertain drift and volatility. The full list of parameters is listed in Table \ref{table:utilityMaxParams}. We present two different case studies.

\begin{table}[tb]
\centering
$$\begin{array}{c} %{|l|r|l|r|}
\hline
r = 0.02, \ T=1, \ \Delta t=0.05   \\ \hline
\alpha = 0.1, \  \gamma = 4, \ \mathfrak{I} = 100, N = 250 \\ \hline
\end{array}$$
\caption{Parameters for the portfolio optimization case study. }
\label{table:utilityMaxParams}
\end{table}

In this setting, we can simplify the computation of $\mathscr{F}_2$ in \eqref{eq:f2} by dimension reduction. Namely,  by monotonicity we
take the 1-dimensional infimum in  the unknown parameters $\mu, \sigma$ just over the boundary $\partial \Theta$ of the parameter space. Let $\rawMu( x^n_{\tk}), \rawSigma( x^n_{\tk})$ denote the parameters achieving the infimum of $\mathscr{F}_2$. Switching to polar coordinates it suffices to find and record the corresponding angle $\varphi^n_{\tk}$ using the map

      %This is achieved by switching to polar coordinates to record the angle $\varphi^n_k$ corresponding to $\rawMu( x^n_k), \rawSigma^2( x^n_k)$ relative to $\procMu^n_k, \procSigma^n_k$:
       \begin{align}
         \rawMu^n_{\tk} &:= \procMu^n_{\tk}+\sqrt{\kappa k^{-1}} \procSigma^n_{\tk} \cos(\varphi^n_{\tk}); \\
         (\rawSigma^n_{\tk})^2 &:= (\procSigma^n_{\tk})^2(1+\sqrt{2 \cdot \kappa k^{-1}}\sin(\varphi^n_{\tk})) \vee 0.
       \end{align}
The minimization defining $\mathscr{F}_2$ (and hence $\varphi^n_{\tk}$) is then done over $\varphi \in [0,2\pi]$ using the default \texttt{fminbnd} algorithm in Matlab with a tolerance threshold of $10^{-6}$. %Overall, this offers a numerical mechanism to generate $\mathscr{F}_2(u, x)$ for arbitrary $u \in \cU, x \in \cX$.
Similarly, we use \texttt{fminbnd} to minimize $\mathscr{F}_2$ over the scalar control domain $\cU = [0,1] = [u_{\min}, u_{\max}]$.

\subsection{Static Investment and Hedging}\label{sec:one-period}
It is instructive to consider the one-period, aka \emph{static} version of the optimal investment problem. In this simplified setting,  the terminal controlled wealth $Y_T$ under $\mathbb{Q}^{\mu,\sigma}$ is given by
\begin{align}
  Y_T^u=Y_{t_0}(1+r(T-t_0)+u(e^{\mu(T-t_0)+\sigma\sqrt{T-t_0}\epsilon_T}-1-r(T-t_0))),
\end{align}
where
$t_0$ is the date that trade happens. Above $\epsilon_T$ is standard Gaussian in $\mathbb{R}$ and due to the scaling of the power utility function, we may take without loss of generality $Y_{t_0} = 1$.
The  static optimal investment problem is
\begin{align}\label{eq:static-inv}
  V(t_0,\procMu,\procSigma)=\sup_{u\in[0,1]}\inf_{(\mu,\sigma)\in\Theta(\procMu,\procSigma,t_0,\kappa)}\mathbb{E}^{\mu,\sigma} \left[\frac{Y_T^{1-\gamma}}{1-\gamma} \right]
\end{align}
with the uncertainty set
\begin{align}
    \Theta(\procMu,\procSigma,t_0,\kappa)=\left\{(\mu,\sigma)\in\mathbb{R}^2:\frac{k_0+1}{\procSigma^2}(\procMu-\mu)^2+\frac{k_0+1}{2\procSigma^4}(\procSigma^2-\sigma^2)^2\leq\kappa\right\}.
  \end{align}
%Note that by default we do not allow short sales or leverage, leading to the constrained control set $\cU = [0,1]$.

Observe that due to the one-period feature, the role of  $\procMu,\procSigma$ is purely to determine the set of $\Q^{\theta}$'s that are being considered. Moreover, without loss of generality we may take the horizon $T-t_0 = 1$ to be fixed, whereby the size of $\Theta(\procMu,\procSigma,t_0,\kappa)$ is controlled by the ratio $\alpha' := \kappa/(k_0+1)$ which is interpreted as the measure of robustness. As $\alpha' \to 0$, $\Theta(\procTh) \to \{ \procTh \}$, and $\alpha'$ is the ``radius'' of adversity.

To solve \eqref{eq:static-inv} it suffices to evaluate the map $(u, \mu, \sigma) \mapsto \mathbb{E}^{\mu,\sigma} \left[\frac{Y_T^{1-\gamma}(u)}{1-\gamma} \right]$ and then find its saddle point (i.e.~the $\sup-\inf$ location) over the constrained domain $u \in \cU$ and $(\mu,\sigma) \in \Theta$. Due to the concavity of the utility preferences, this map is increasing in $\mu$ and decreasing in $\sigma$. Consequently, the infimum over the robustified beliefs will always be achieved in the NW quadrant, so that $\check\mu(\procMu,\procSigma) \le \procMu$ and $\check\sigma(\procMu,\procSigma) \ge \procSigma$. Therefore we can again switch to polar coordinates reducing to a 2-dimensional function $F(u,\varphi)$ of the control $u$ and the angle $\varphi \in [\pi/2, \pi]$. For example, $\varphi = \pi$ means taking the minimum possible drift.

\begin{figure}[ht]
  \centering
    \includegraphics[width=0.45\textwidth]{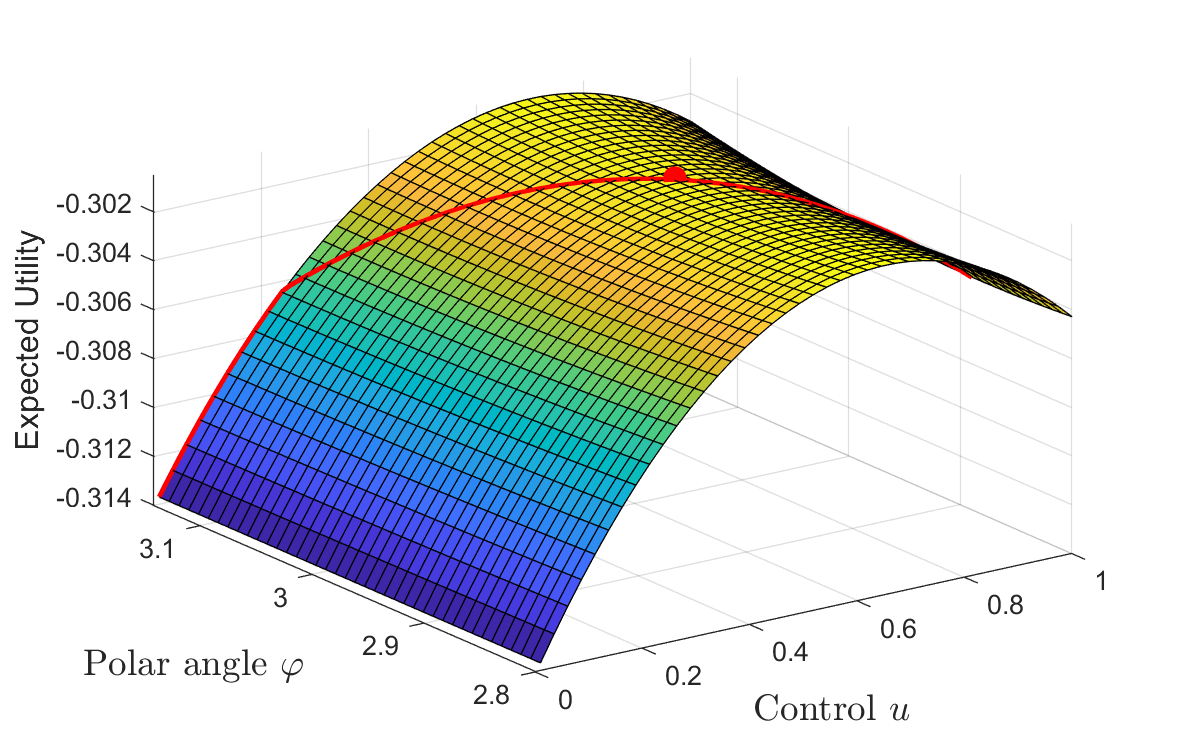}  \includegraphics[width=0.45\textwidth,trim=0in 0in 0in 0.5in]{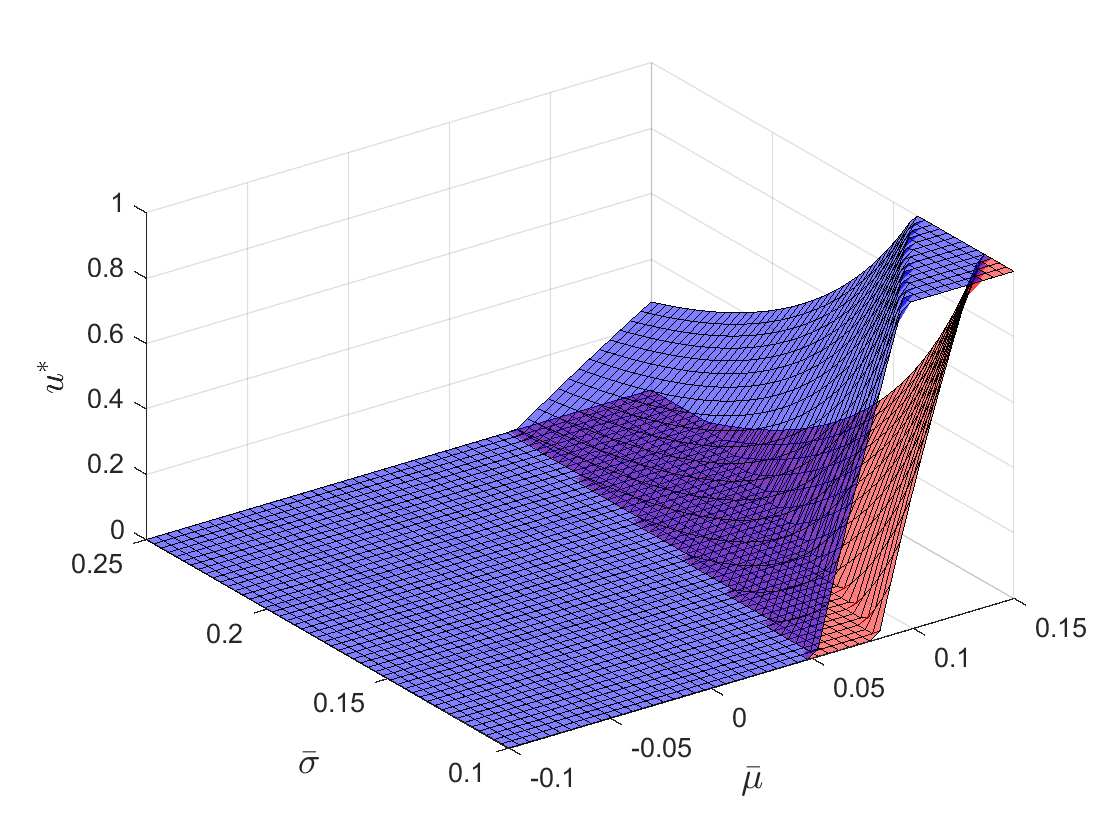}
 \caption{ {\em Left:} $F(u,\varphi)$ for the static portfolio optimization case study for a fixed $\procMu,\procSigma$.  The red line denotes the profile minimizers $u \mapsto \check\varphi(u)$. The red dot denotes the saddle point $(u^*, \varphi^*)$ identified via $u^* = \arg\max_{u \in [0,1]} F(u, \check\varphi(u))$. \emph{Right:} $u^* \in [0,1]$ as a function of $\procMu,\procSigma$ for two different robustness levels $\kappa=4.61$ (red), $\kappa = 1.39$ (blue) and $k_0 = 10$. \label{fig:saddle-point}}
\end{figure}

The left panel in Figure \ref{fig:saddle-point} shows $F(u,\varphi)$ for a representative pair of initial beliefs $\procMu,\procSigma$ in the range $u \in [0,1], \varphi \in [2.8,\pi]$. For a fixed investment level $u>0$, $F(u,\cdot)$  is a convex function with a minimizer $\check\varphi(u)$. However note that when $u =0$, $Y_T$ is independent of the $S$-dynamics, i.e.~$F(0,\cdot)$ is constant in $\varphi$. The respective infimum is then undefined and there is no $\check\varphi(0)$. The saddle-point is the maximum $\arg\max_u F(u, \check\varphi(u))$ and is indicated by the red dot on the left panel, where $u \mapsto F(u, \check\varphi(u))$ is visualized by the red line. In this particular case, $u^*(\procMu,\procSigma) \simeq 0.669$ and $\varphi^*(\procMu,\procSigma) \simeq 2.998$.

The right panel of  Figure \ref{fig:saddle-point} displays $u^*$ as a function of the initial beliefs $\procMu,\procSigma$ for two different robustness levels. As $\alpha$ decreases, $\Theta$ gets larger and therefore the worst-case (inner infimum) becomes worse. Consequently, the controller acts more conservatively, i.e.~$u^*$ increases in $\alpha$. Otherwise, $u^*$ has the familiar shape of increasing in $\procMu$ (asset returns being more favorable) and decreasing in $\procSigma$ (asset risk rising). Note that when returns are not sufficiently high and/or volatility is too large, the optimal action is $u^*=0$, i.e.~to invest only in the risk-free bond. This happens in particular when $\Theta(\procMu,\procSigma,t_0,\kappa) \bigcap \R^2_- \neq \emptyset$ intersects with the negative half-plane, i.e.~the investor believes that negative returns $\mu < 0$ are possible.

\begin{figure}[ht]
  \centering
    \includegraphics[width=0.45\textwidth,height=2in]{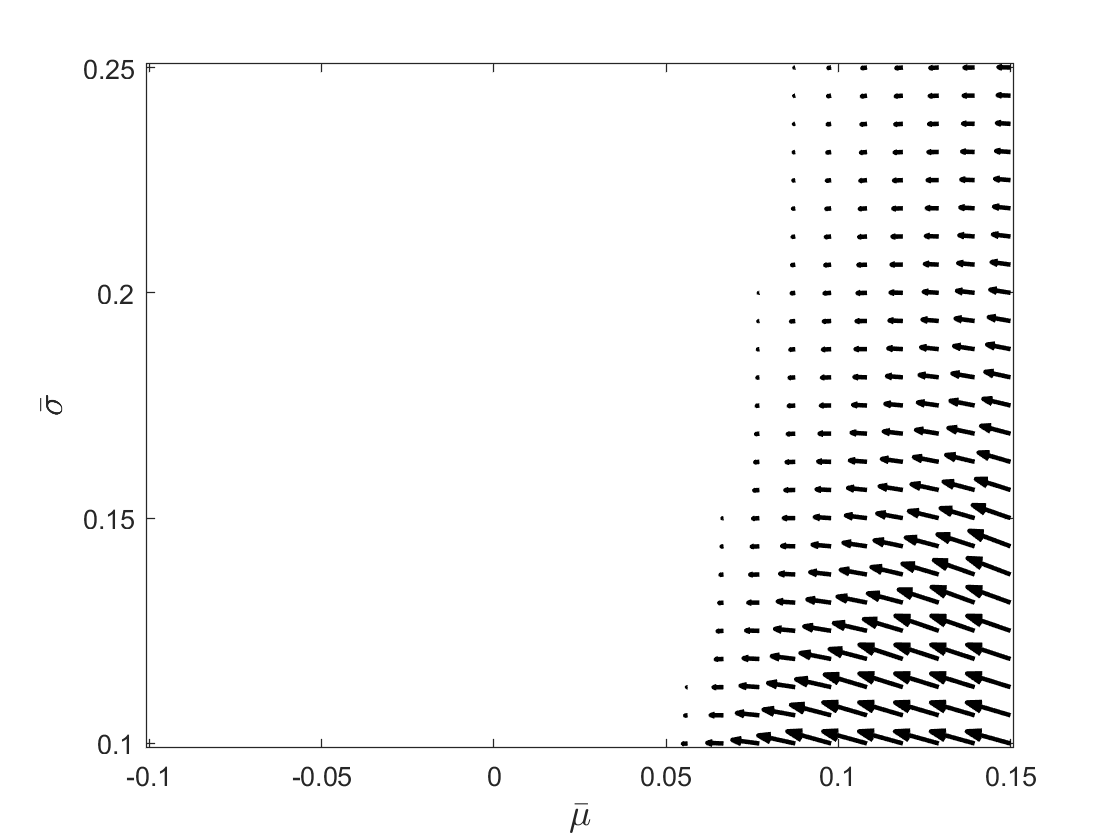}  \includegraphics[width=0.45\textwidth,trim=0in 0.2in 0in 2in]{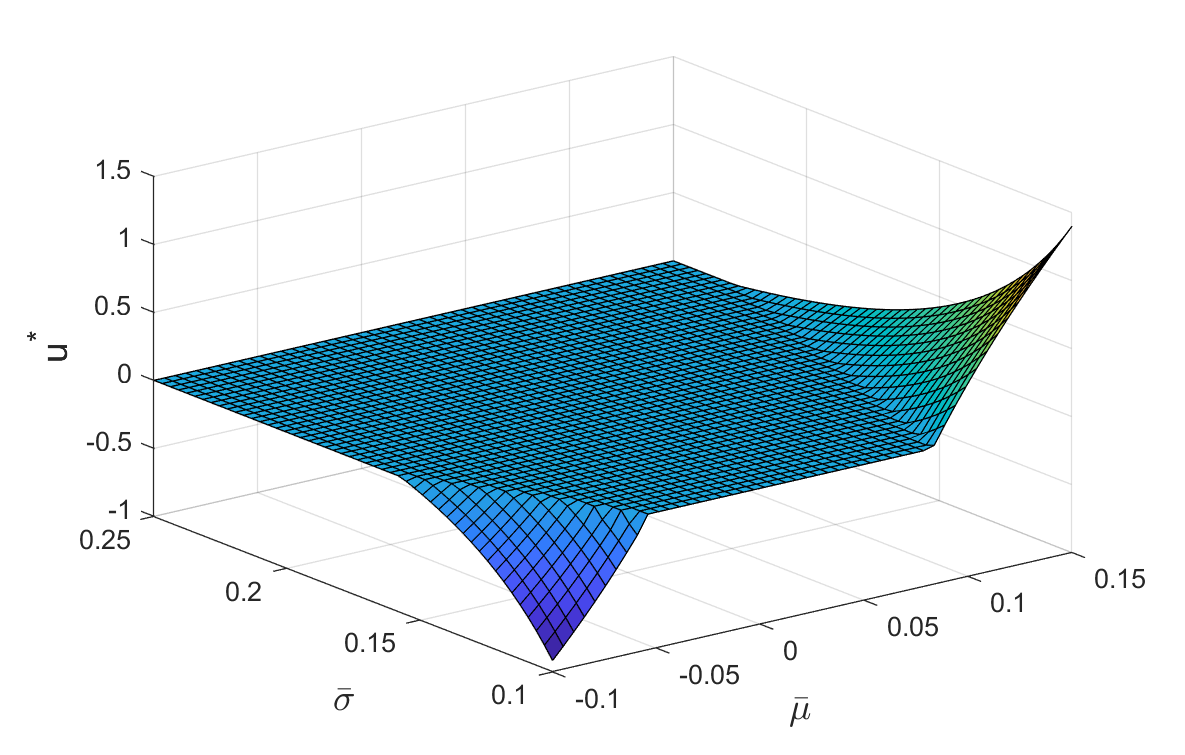}
 \caption{{\em Left panel:} quiver plot showing the relative position of $(\check\mu, \check\sigma)$ relative to $(\procMu, \procSigma)$. The angle of the quivers is the 
robustified belief angle $\varphi^*$ and the length is proportional to $u^*$; no quiver is shown  when $u^* = 0$. \emph{Right:} optimal investment fraction $u^*$ as a function of $\procMu,\procSigma$ using unconstrained optimization over $\mathbb{R}$. \label{fig:static-phi}}
\end{figure}

The left panel of Figure~\ref{fig:static-phi} visualizes the structure of the robustified beliefs $(\check\mu, \check\sigma)$ as a function of $(\procMu, \procSigma)$. Recall that $\varphi^* = \pi$ (quivers pointing to the West)  corresponds to the worst-case making $\check\mu$ as negative as feasible. This is the robustified belief in the middle and top of the state space. In the bottom-right (very high $\procMu$ and very low $\procSigma$), the worst-case trades off reduced returns against increased risk. Finally, in the left half of the plot we have $u^* = 0$, so that  there is no $\varphi^*$ (no quivers displayed) since $F(0,\cdot)$ is constant. In that case the consistent interpretation is to keep  $\check\mu = \procMu, \check\sigma = \procSigma$ unchanged. %Overall, we confirm a strong negative relationship between $\varphi^*$ and $u^*$: when the investor is taking large positions in the risky asset she is increasingly concerned about high volatility rather than just low returns. Note that when

%We are interested in studying the optimizer map $(\hat{\mu},\hat{\sigma},t_0,\kappa)\mapsto(\mu^*(\hat{\mu},\hat{\sigma},t_0,\kappa),\sigma^*(\hat{\mu},\hat{\sigma},t_0,\kappa))$ that is linked to the saddle-point of $F(\cdot, \cdot)$.

%\blu{?? take $\Theta$ as a general set valued function other than the $(1-\alpha)$-confidence region, and study the dependence of $(\mu^*,\sigma^*)$ on the position, shape and size of $\Theta$.} \textbf{Effect of $\Theta$: }
%Change the size of $\Theta$ via $\alpha$ and the shape (rectangle vs ellipse). The worst-case is generally to lower $\mu$ and to increase $\sigma$. Hence the minimizer will be trivial (at the extremal point or corner) for a rectangular $\Theta$.

Finally, the right panel of Figure~\ref{fig:static-phi} addresses the question of trading constraints. Typically we find that $u \mapsto F(u,\check\varphi(u))$ is convex, i.e.~has a unique global maximum. In that case we may find the unconstrained global maximum $u^*_{unc}$ and then constraints of the form $u \in [u_{\min}, u_{\max}]$ translate into $u^*_{con} = u_{\min} \vee u^*_{unc} \vee u_{\max}$. The plot reveals several useful insights. First, we observe the ``shrinkage to zero'' effect: robust control implies that a risky position will be assumed only if there is a clear understanding regarding $\mu$. Therefore, if $| \procMu - r|$ is small, the robust solution is to take $u^* = 0$. In that sense, $u^*$ shrinks the non-robust control $u^*_{adpt}$ to zero. Because the size of $\Theta(\tk,)$ is influenced by $\hat{\sigma}$ the amount of this shrinkage grows in $\hat{\sigma}$ as observed in the plot. Second, we observe that for a fixed $\procSigma$ $u^*$ is essentially piecewise linear in $\procMu$, capturing the original Merton intuition that $u^*$ is linear in asset returns. Third, we note that leverage and short-selling remain feasible with adaptive robust control provided that the adversarial set is far enough to one side of the half-plane. Thus, we observe $u^* < 0$ in the bottom-left corner (strongly negative Sharpe ratio) and $u^* \gg 1$ in the bottom-right corner (very high Sharpe ratio). Fourth, we stress that the plot is (roughly) symmetric about $\procMu = r$ which is the risk-free alternative.

%Figure \ref{fig:saddle-point} also demonstrates that frequently the constrains of $u^* \in [0,1]$ are binding, in other words on the interval and hence achieves its maximum at one of the two boundaries. This structure suggests the possibility of doing unconstrained optimization and then projecting back into the interval as a way of numerically smoothing this effect. \blu{??? Consider the possibility of short-selling -- how does this modify $u^*$? Should not much, as need to be really convinced that $\mu < 0$ to short-sell. }

%\blu{Add a figure for $\varphi^*(\hat{\mu},\hat{\sigma})$. $\varphi^*$ is undefined when $u=0$}

\subsection{Multi-period Optimal Investment}

\begin{figure}[ht]
  \centering
    \includegraphics[width=0.45\textwidth]{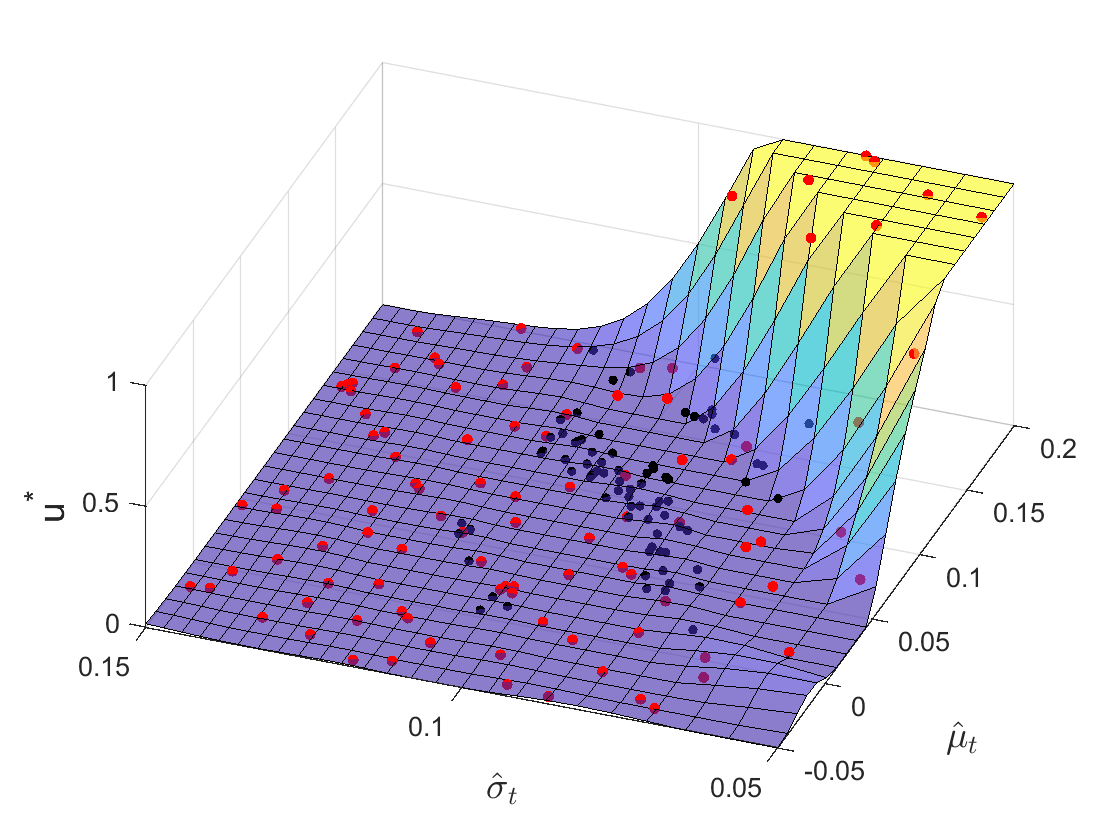}
    \includegraphics[width=0.45\textwidth]{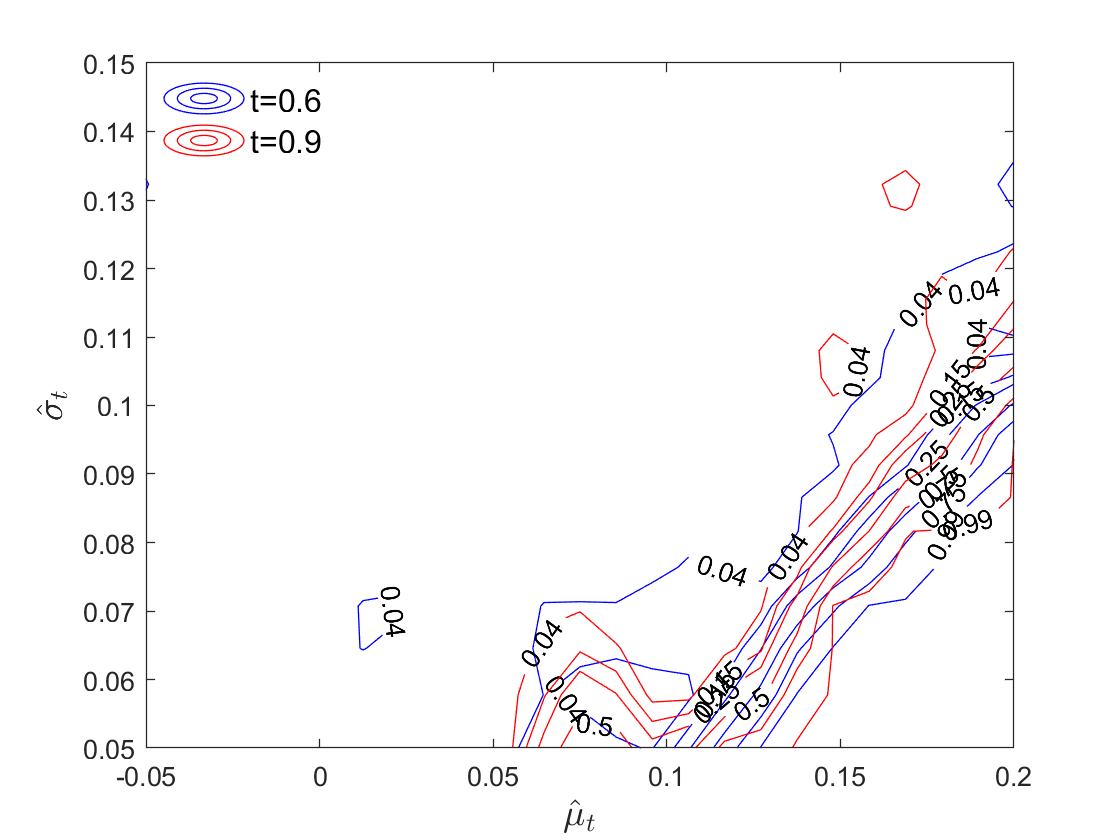}
\caption{Left: Optimal investment strategy from the GP surrogate, $\modelU({\tk},x)$ at ${\tk}=t_0+15\Delta t$. Right: time-dependency of $\modelU({\tk},\cdot)$: comparing $k=15$ (cyan) and $k=10$ (blue). \label{fig:opt-a}}
\end{figure}

We now return to the original multi-period setting with $K=20$ time-steps, i.e.~$t_k = t_0 + k\Delta t$ and $T= t_K$. This setup is very similar to the previous section, except that $\modelU(\tk,x)$ is now time-dependent.   Figure \ref{fig:opt-a} illustrates this dependence by showing a contour-plot of $(\procMu_\cdot, \procSigma_\cdot) \mapsto \modelU(\cdot,\procMu, \procSigma)$ for two different time steps $t_k$. Recall that at early epochs $\Theta_{k}$ is larger so ceteris paribus the investor is more conservative and  $\modelU$ tends to increase in $\tk$. This is indeed observed in the figure, matching the intuition of ``learning-as-you-go''. We note that there are competing effects, in particular due to the time-dependence in the dynamics of $\procTh_{\tk}$ (learning slows over time) and the time-dependent effective risk-aversion (the value function becomes less concave as $T-{\tk}$ increases).

\textbf{Structure of Optimal Control: }
The left panel of Figure \ref{fig:opt-a} shows the optimal feedback control surface $(\procMu, \procSigma) \mapsto \modelU({\tk},x)$. We observe that $\modelU({\tk},x) = 0$ when $\procMu$ is low and monotonically increases as the posterior mean asset return $\procMu_{\tk}$ rises. Eventually for high enough $\procMu_{\tk}$, $\modelU({\tk},x) = 1$. This structure suggests that we can concentrate the statistical modeling efforts in the intermediate region $R$ where $\modelU({\tk},x) \in (0,1)$ since otherwise the corresponding feedback control surface is completely flat and therefore easy to interpolate. The latter fact also implies that even extrapolation is feasible since we just need to ensure that the surrogate is set such that the respective asymptotes in $\procMu$ are $0$ to the left and $+1$ to the right.

\textbf{Structure of Robustified Beliefs: } We recall that $\hat{\varphi}(\tk, x)$ is only well-defined when $\modelU(\tk,x) \neq 0$, since otherwise the investor is not exposed to returns uncertainty. Therefore, $\hat{\varphi}(\tk,x)$ only exists for $\procMu_{\tk}$ high enough and $\procSigma_{\tk}$ low enough (the SE corner of the state space). In the latter region, we find that the principal effect is to decrease the drift rather than to increase volatility, so that $\hat{\varphi} \simeq \pi$. In the multi-period context, low $\procSigma_{\tk}$ would also hurt the investor, so $\hat{\varPhi} > \pi$ will also occur at the early time steps.

%For the optimization over $\Theta(t, \procMu, \procSigma)$, we only optimize over its boundary; since the latter is an ellipse in the base setup, we implement dimension reduction by switching to polar coordinates, mapping $\partial \Theta$ into an \emph{angle} $\varphi \in [0,2\pi]$. For example, $\varphi = \pi$ means taking the minimum possible drift.
%Figure \ref{fig:opt-angle} shows estimated $\rawPhi$ which we see is somewhat concentrated around $\pi$, i.e.~the principal effect is to decrease the drift, rather than to increase volatility. We also note that in the multi-period context, low $\sigma$ would also hurt the investor so frequently $\rawPhi > \pi$. \blu{Discuss dependence on $t$}
%
%\begin{figure}[ht]
%  \centering
%    \includegraphics[height=2in]{RegressionPictures/angle_quiverPlot.png}
%\caption{Estimated angle $\rawPhi$ for the robustified model parameters $\rawTh$. Because the size of $\Theta(t,\procTh_t)$ depends on $\procTh_t$, the difference $\rawTh-\procTh_t$ is state-dependent and indicated by the quiver plot.  \label{fig:opt-angle}}
%\end{figure}

\subsection{Distribution of Terminal Wealth}
To obtain an out-of-sample performance metric of the computed investment strategy $\modelU$ we apply Algorithm \ref{algo:fwd} in Section~\ref{sec:algo}. For the test measure $\Q^{test}$ we use $\Q^{\Theta_{0}}$, i.e.~static Knightian uncertainty where $\theta^{*,n}$ is drawn from a known prior $\Theta_{0}$ at $t=t_0$ and then kept constant through time. To be consistent, we utilize the initial adversarial set $\Theta_{0} =\Theta_{0}(\alpha; \procTh_{t_0})$. Specifically, for the forward simulations we sample independently and uniformly $\theta^{*,n}$ for $n=1,\ldots, N'$. %This nested procedure with a total of $M' = N_{in} \times N_{out}$ simulations then allows us to evaluate the conditional $Y^{\vec{u}}_T | \theta^{*,(n)}$ and the overall $\E^{\Q^{\Theta_0}}[ U(Y^{\vec{u}}_T) ]$. It also captures the learning aspect as $\procTh_t \simeq \theta^{*,(n)}$ for each inner batch.
For the discussion below, we took $\mu^{*,n} \sim \mathcal{N}(0.15, 0.02^2), \sigma^{*,n} = 0.1$ and $\procTh_{t_0} = (0.1, 0.08)$ so the investor starts out under-estimating both the mean returns and the return volatility.  %(0.15-0.02)/( 4*0.01) = 0.13/0.04 > 1 if knew; start out with (0.1-0.02)/(4*0.08^2) = 1/(4*0.08)

We recall that the classical Merton solution with fixed $(\mu, \sigma)$ is time stationary and given by
$$
u^*(t, x) = u^*(\mu, \sigma) = \frac{\mu-r}{\gamma \sigma^2}.
$$
Using the above we can compare our adaptive robust solution to the following alternative strategies:
\begin{itemize}
  \item Merton strategy based on $u^*(\modelTh({t_0}, \procTh_{t_0}))$ which is the static robust formulation;

  \item Merton strategy based on $u^*(\procTh_{\tk})$ which is the adaptive formulation;

  \item Strategy equivalent to our adaptive robust algorithm with $\alpha =1$ ($\Theta_{k}(\procTh_{\tk}) = \{ \procTh_{\tk}\}$) which dispenses with robustness.
\end{itemize}
For the static robust approach, the worst case of asset dynamics corresponds to a low Sharpe ratio. In particular, for our parameter values, since negative return rate cannot be ruled out $\Theta_0 \bigcup (-\infty,0) \times \R_+ \neq \emptyset$, the robust $u^{(SR,*)}(t,x) = 0$ for all $t\in\cT\setminus\{T\}$, and no investment would be ever undertaken. This point clearly illustrates why learning is desirable, as it allows the investor to overcome initial pessimism (or worry about model risk) through gradual collection of information. If enough favorable information is learned then she will eventually become optimistic enough to invest. At the other extreme, the adaptive case leads to over-confidence and hence over-investment. In particular, the adaptive method generates much larger variance, i.e.~risk, in $Y_T$ compared to adaptive robust, as the investor trusts her myopic beliefs too much.

%We then compare the above distribution to one based on an adaptive strategy (that applies $u( \procTh_t)$ rather than $u( \rawTh( \procTh_t) )$ ) and the Merton strategy based on $u( \theta_0)$.
%\blu{Also show the impact of $\alpha$. Show a histogram of $Y_T$ + 2 tables}.

\begin{figure}[ht]
  \centering
    \includegraphics[width=0.45\textwidth]{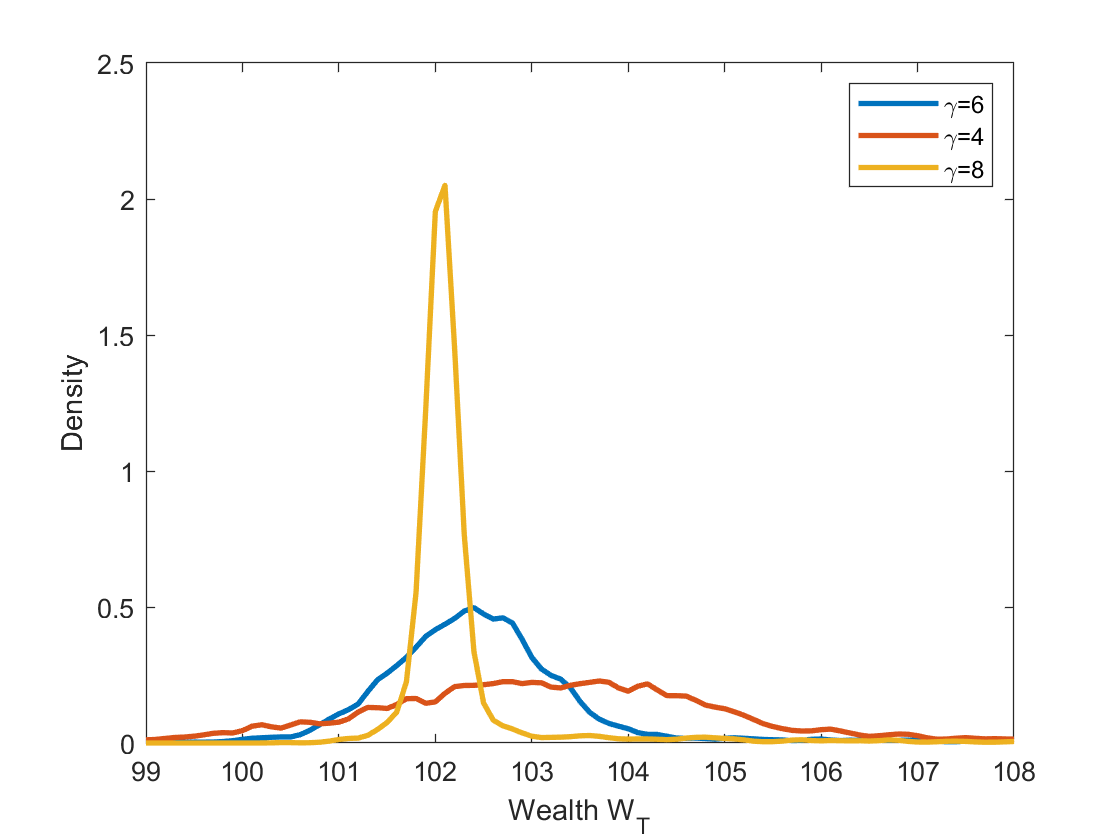}
    \includegraphics[width=0.45\textwidth]{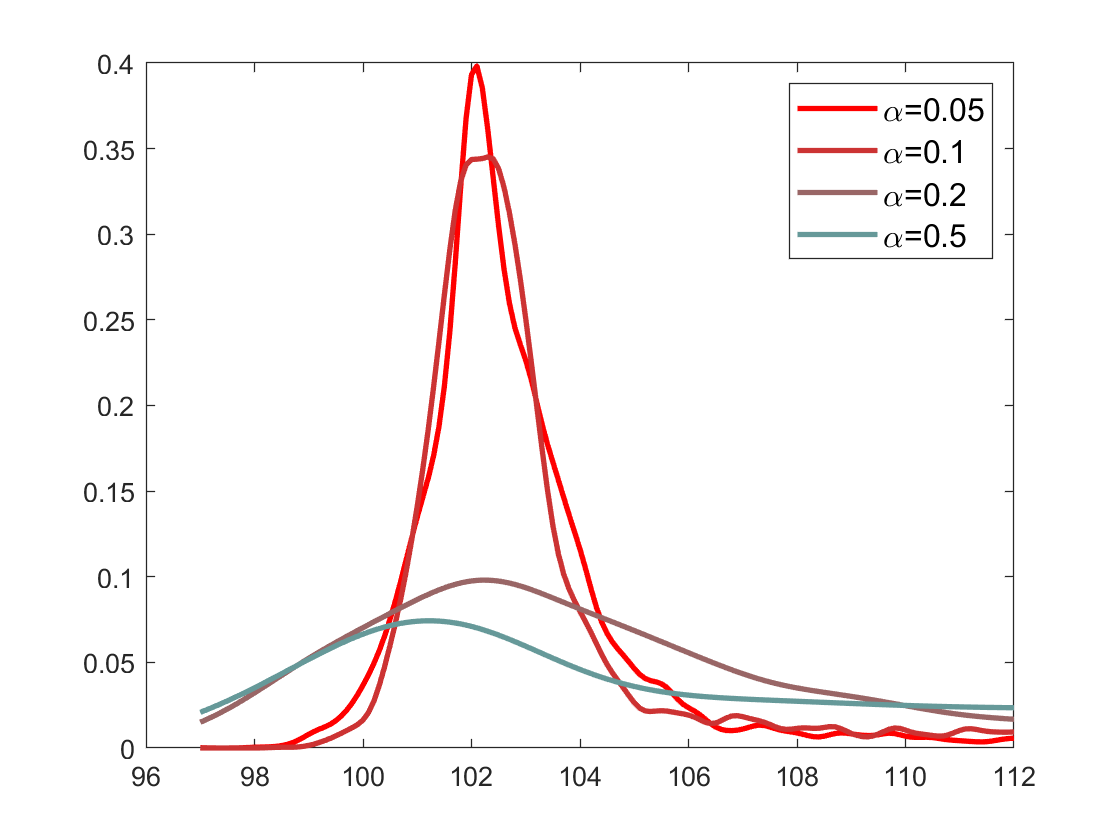}
\caption{Left: Distribution of $W_T$ for different risk aversion parameters $\gamma$. Right: Distribution of $W_T$ for different robustness parameters $\alpha$. \label{fig:opt-w}}
\end{figure}

Figure \ref{fig:opt-w} shows the distribution of terminal wealth $W_T$ as we vary the risk aversion parameter $\gamma$ and the robustness parameter $\alpha$. As expected, less risk averse investors will have higher $u^*$ and therefore larger $\E[W_T]$ accompanied by much larger $\Var(W_T)$. Similarly, agents that are less conservative (and therefore have smaller $\Theta_{k}(\procTh_{\tk})$) will tend to again invest more. However, being too ``trusting'' and not guarding against model mis-specification will eventually
hurt wealth accumulation (compare $\alpha=0.5$ to $\alpha=0.2$, so that $\E^{\Q^{\Theta_0}}[ W_T]$ has a hump shape as a function of $\alpha$).

%\subsubsection{Adaptive Control}
%In the adaptive control formulation, there is myopic learning but no robustness. In other words, we preserve the dynamics \eqref{eq:T-optInv} but take $\Theta(t,\procMu, \procSigma) = \{ (\procMu, \procSigma) \}$ as a singleton, eliminating any model uncertainty. Practically this means substituting $\procMu, \procSigma$ for $(\mu,\sigma)$ in \eqref{eq:tilde-V} and dropping the inner infimum. We can then treat $\procMu, \procSigma$ as optimization parameters and only need to solve the \emph{parametric}, time-stationary optimization problem
%\begin{align}
%  \check{V}( \procMu, \procSigma ) = \sup_{u \in [0,1]} \E[ (1+r+u(\procMu + \procSigma \eps )^{1-\gamma} ) ]
%\end{align}
%which then leads to $\tilde{V}(t,\procMu, \procSigma) = \check{V}(\procMu, \procSigma) \E[ \tilde{V}(t+1,\procMu_{t+1}, \procSigma_{t+1}) | \procMu, \procSigma] $.
%
%In other words, the adaptive solution consists of two phases. During the optimization phase, we solve the Bellman equations with $\rawMu$ replaced by $\mu$, and for all $\mu\in \Theta$. The optimal selector is denoted by $u_t^\mu, \ t\in\cT'$, which is in fact independent of $t$. In the second adaptation phase, for every $t\in \{0,1,2,\ldots,T-1\},$ we compute the pointwise estimate $\procMu_t$ of $\theta ^*$, and apply the certainty equivalent control  $u_t=u^{\procMu_t}_t$. For more details see, for instance, \cite{KumarVarayia2015Book}. %\cite{ChenGuo1991-Book}.
%

\section{Adaptive Robust Optimal Hedging}\label{sec:hedging}

It is of financial and theoretical importance to study optimal hedging problem in an incomplete market setup under model uncertainty. The uncertainty issue was addressed, for example, in the paper \cite{monoyios07} where the author used Bayesian estimators in place of the corresponding parameters. We will compare our approach to a similar method that uses maximum likelihood estimators instead. Another standard method widely used among banks when dealing with model risk is the so-called conservative or worst-case approach, and it is in principle the same as the static robust approach.

We denote by  $S=\{{S}_t,\, t\in\cT\}$ the  market price process of the stock  and we postulate that  the increments of  $S$ are
log-normal with parameters $\mu,\sigma$:
%Equivalently, the log-asset process follows
\begin{align}\label{eq:S-dyn}
S_{\tkp} = S_{\tk} \exp \left( \mu \Delta t + \sigma \sqrt{\Delta t} \eps_{\tkp} \right), \qquad \eps_{\tkp} \sim \mathcal{N}(0,1).
%Y_{t+1} = Y_t + \mu + \sigma \eps_t, %\qquad Y_t := \log S_t.
\end{align}
Dynamic trading strategies are based on a discount bond and the stock $S$, and are represented by a predictable discrete-time process  $\{u_t,\ t\in \cT'\}$ identified with the number of shares invested in the stock from ${\tk}$ until $\tkp$. The bond investment is then determined from the self-financing constraint.  For an admissible $(u_t)$ we define the corresponding wealth process $W$ as
\begin{align}\label{eq:W-dyn}
   {W}_{\tkp}
   = W_{\tk} + u_{\tk}  ( S_{\tkp}-  S_{\tk})=W_{\tk}+ u_{\tk} S_{\tk}(e^{\mu \Delta t+\sigma \sqrt{\Delta t}\epsilon_{\tkp}}-1), \quad   W_{t_0}=w,
\end{align}
where $w$ is the initial endowment, assumed to be fixed and deterministic.

We are interested in dynamic hedging, using an admissible trading strategy, of a European style option written on the stock $S$, maturing at time $T$ and with the payoff $\Phi(S_T)$, where $\Phi$ is a measurable function. To address nonlinear hedging, we define the objective as optimization of the hedging \emph{error} $H = \Phi(S_T) - W_T$ through \cite{RTV2002}
\begin{equation}\label{eq:riskmeasure1}
\mathbb{E} \left[  \ell(\Phi( S_T) -  W_T)   \right ] \mapsto \min!
\end{equation}
where $\ell(\cdot)$ is a positive measurable \emph{loss function} such that $\ell(0)=0$. This setup offers a loss-based criterion that is applicable both in complete and incomplete markets. The main loss function we will consider is of the form
\begin{equation}\label{eq:riskmeasure2}
\ell(h)= h^++\lambda h^-,
\end{equation}
where the constant $\lambda\geq0$ reflects the relative importance of super- and sub-hedging. When $\lambda =0$ the investor is penalized only for not having enough to cover her liability, while $\lambda = 1$ leads to $\ell(h) = |h|$ and is the $L^1$-equivalent of the classical quadratic hedging criterion.

We consider  adaptive robust hedging in the context of uncertain asset drift and volatility, leading to a four-dimensional state: $x=  (W,S,\procMu, \procSigma)$ which under $\Q^{\mu,\sigma}$ has the dynamics prescribed by \eqref{eq:S-dyn}-\eqref{eq:W-dyn}-\eqref{eq:mu-dyn}-\eqref{eq:sigma-dyn}.
The Bellman equation becomes
\begin{align}\notag
  V(T,x) &= \ell(\Phi(S)-W),\\ \label{eq:hedge-bellman}
  V({\tk},x) &= \inf_{u \in \cU} \sup_{ \theta=(\mu,\sigma) \in \Theta({\tk},\procMu, \procSigma)} \E \left[ V(\tkp, \bT({\tk},x,u,\theta,\eps_{\tkp}) ) \right].
\end{align}
We may interpret $V$ as the robust expected loss, so that $V({\tk},x) \geq 0$ and lower $V$ implies better ability to $\ell$-hedge. For example, in the one-sided example with $\ell(h) = h^+$, $V({\tk},\cdot)$ will be monotone decreasing in both $S$ and $W$.

This setup is conceptually similar to the previous section but
presents a harder computational challenge due to its higher dimensionality and the non-smooth relationships between $V({\tk},\cdot)$ and the state variables. Moreover, it features the fully endogenous wealth process $W$ which cannot be reasonably ``forward-simulated'' without assuming $u_{\tk}$. To overcome this issue we adopt the Control Randomization approach of selecting a $\mathbb{Q}^{0}$ measure to build the designs $\cD_{k}$.

For illustration, we henceforth consider a European style Call option written on the risky asset $S$ with strike ${\cal K}$,  $\Phi(S)=(S-\mathcal{K})^+$. The investor is short (sold) the Call and aims to minimize her hedging error starting with initial endowment $W$.

\subsection{Implementation Details}

\begin{table}[tb]
\centering
$$\begin{array}{c} %{|l|r|l|r|}
\hline
r = 0, \ T=1, \ \Delta t=0.1, \ k_0=150 \\ \hline
\alpha = 0.1, \  \kappa = 4.61, \ \mathfrak{I} = 40, \ N = 250 \\ \hline
\bar{\mu}_{t_0}=0.12, \ \bar{\sigma}_{t_0}=0.4, \ S_0=100, \ \mathcal{K}=100 \\ \hline
\end{array}$$
\caption{Parameters for the optimal hedging case study. }
\label{table:HedgingParams}
\end{table}

% As mentioned earlier, due to the higher dimensionality and lack of dimension-reduced Bellman equations, the hedging problem poses a more difficult challenge in applying the adaptive robust methodology. Therefore, it is crucial to choose a proper simulation design $\cD_k$ in order to make the problem tractable from computational point of view.

For the hedging problem, $x$ is four-dimensional, and there is a high correlation between $S_{\tk}$ and $\bar{\mu}_{\tk}$, so that the design $\mathcal{D}_{\tk}$ has a certain 3-d structure. Moreover, necessarily there will be a dependence between $S_{\tk}$ and wealth $W_{\tk}$.
To construct $\cD_k$ we use a three-step procedure. First, we generate ``pilot'' forward paths of $X$ under $\mathbb{Q}^{0}$. The pilot paths are based on randomized $\procMu_{t_0}^{1:N_0}$, $\procSigma_{t_0}^{1:N_0}$, as well as initial stock prices $S_0^{1:N_0}$, where $N_0=250$. Specifically, given $\mu_0, \sigma_0, K$ we sample initial beliefs $\procMu_{t_0}^{1:N_0}$, $\procSigma_{t_0}^{1:N_0}$ uniformly from $[0.5{\mu}_0,1.5{\mu}_0]$ and $[0.6{\sigma}_0,1.3{\sigma}_0]$, respectively, as well as $S_{t_0}^{1:N_0} \in [0.5\mathcal{K}, 2\mathcal{K}]$ around the strike price  $\mathcal{K}$.
Second, we construct a set of augmented pilot sites by adding a few more $(\procMu, \procSigma)$'s  based on the edges of adversarial beliefs $(\mu,\sigma)\in\partial\Theta({\tk},\bar{\mu}^n_{\tk},\procSigma^n_{\tk})$. This step helps us reduce extrapolation when evaluating $V(\tkp, \bT({\tk},x,u,\theta,\eps^{(i)}))$ during Gaussian quadrature.
Third, we use a space-filling strategy to fill in ``holes'', employing a QMC sequence-based design over the convex hull of the augmented pilot sites at each $t_k$.
%
%
%For easing the burden of the high dimensionality, our main idea is to use the CR approach to generate forward paths of the process $(S,\bar{\mu},\procSigma)$ under some chosen probability measure $\mathcal{Q}^0$. This allows us to avoid the states that are unlikely to be visited. But, several issues needs to be addressed regarding the sequential design.
%
%Due to the high correlation between $S$ and $\bar{\mu}$, a naive forward simulation of the augmented state process will result in a sequence of 3-dimensional hyperplane in the 4-dimensional state space. Such phenomenon then prevents the application of space-filling and can lead to unfavorable numerical results.
%As a solution, as those are the really interesting ones. By doing so we are able to better utilize the space-filling technique to achieve a good approximation of the continuous state space.
%
This is done by ``pulling forward'' the design $\cD_{k}$ to the next time-step $t_{\tkp}$ using the map $\mathbf{T}({\tk},x^n_{\tk},u,(\mu,\sigma),\epsilon^{(i)}_{\tkp}))$ with the empirically sampled $\epsilon$'s and then randomly replacing 100 pilot sites by sites from the respective QMC (Sobol) sequence in $\R^3$. This yields the experimental design $(S^{1:N}_{\tk},\bar{\mu}^{1:N}_{\tk},\bar{\sigma}^{1:N}_{\tk})$, $N=250$.

%For the space-filling part, we approximate the hypersurface of $(S_{k+1},\bar{\mu}_{k+1},\procSigma_{k+1})$ by using 150 sites from the convex hull generated by
%\begin{align*}
%  S^{(i)}_{k+1}&=S^{(i)}_ke^{\mu+\sigma\epsilon^{(i)}},\\
%  \mu^{(i)}&=\frac{t+t_0}{t+t_0+1}\bar{\mu}_k^n+\frac{1}{t+t_0+1}(\mu+\sigma\epsilon^{(i)}),\\
%  \sigma^{(i)}&=\sqrt{\frac{t+t_0}{t+t_0+1}(\procSigma^{n}_k)^2+\frac{t+t_0}{(t+t_0+1)^2}(\bar{\mu}^n_k-\mu-\sigma\epsilon^{(i)})^2},\\
%  i&=1,\ldots,\mathfrak{I},\ n=1,\ldots,N_k,\ (\mu,\sigma)\in\partial \Theta(t,\bar{\mu}_k^n,\procSigma_k^n),
%\end{align*}
%as well as points on the forward simulated paths. To improve such approximation, we then select 100 sites through a QMC space-filling technique to remove the holes. %\textcolor[rgb]{1.00,0.00,0.00}{see picture}

Finally, to specify $W^{1:N}_{\tk}$ we compute the Black Scholes prices $P^{BS}({\tk},S^{n}_{\tk}; \procTh^n_{\tk})$ and independently sample $W^n_{\tk} \sim Unif( 0.5 P({\tk},S^n_{\tk}), 1.5 P^{BS}({\tk},S^n_{\tk}))$. This leads to a ``tube'' in the $(S,W)$ plane that contains the Black Scholes prices, with the idea that these are the likely wealth levels (assuming that $W_{t_0} \approx P^{BS}(t_0,S_{t_0}; \procTh_{t_0})$) to be visited.

%by perturbing from the time $t$ Black-Scholes price of the call option $P(S_t^n)$ by no more than 50 percent.
%where $\mathrm{BS}(S_t^n)$ is the Black-Scholes value of the option at time $t$ for stock price $S_t^n$.
%This yields an economical and efficient discretization for the space of $W$.
%\textcolor[rgb]{1.00,0.00,0.00}{maybe try to use one picture to show such point}
%Design $\cD_k$, $k=1,\ldots,K$, is then a ``tube'' that wraps around states that are likely to be visited under the measure $\mathcal{Q}^0$.

%\begin{figure}[ht]
%  \centering
%    \includegraphics[width=0.45\textwidth]{RegressionPictures/simDesignHedging.png}
%    \includegraphics[width=0.45\textwidth]{RegressionPictures/simDesignV.png}
% \caption{Left: constructed simulation design $\cD^{\procMu,\procSigma,y} = \{ \procMu^{1:N}, \procSigma^{1:N},y^{1:N} \}$ for the optimal hedging case study. There are a total of 250 sites: density based (red); QMC based (blue). Right: two designs $\cD^{\procMu,v}=\{ \procMu^{1:N},v^{1:N} \}$: forward simulation based (red); perturbation of Black-Scholes price based (blue). }
%\end{figure}

%\subsection{Algorithm for Adaptive Robust Hedging}

The algorithm described in Section~\ref{sec:algo} is tailored for applying the adaptive robust framework to the optimal hedging problem.
For the hedging problem monotonicity of $\hat{E}[\hat{V}(\tkp,\mathbf{T}({\tk},x^n_{\tk},u,(\mu,\sigma),\epsilon_{\tkp}))]$ in terms of $\mu$ or $\sigma$ no longer holds true.
Hence, the supremum is not necessarily attained on the boundary $\partial\Theta({\tk},\bar{\mu}^n_{\tk},\procSigma^n_{\tk})$. Nevertheless, we continue to take advantage of the elliptical nature of $\Theta({\tk},\cdot)$ working in the respective polar coordinates $(\varphi, \rho)$ representing the angle and distance from $(\procMu,\procSigma)$. We then discretize  $\Theta({\tk},\bar{\mu}^n_{\tk},\procSigma^n_{\tk})$ uniformly in terms of $\varphi \in [0, 2\pi]$ and $\rho \in (0, \kappa)$ and fixing $u \in \cU$ carry out a direct maximization over the discrete set $\Theta'$ of adversarial beliefs. This finally yields the robustified parameters $(\rawMu,\rawSigma)$:
\begin{align*}
    \rawMu^n_{\tk} &= \procMu^n_{\tk}+\sqrt{\rawRho^n_{\tk} \cdot (\procSigma^n_{\tk})^2/(k+k_0)} \cos(\rawPhi^n_{\tk}) \\
    \rawSigma^n_{\tk} &= \procSigma^n_{\tk} \sqrt{1+\sqrt{2 \cdot \rawRho^n_{\tk}/(k+k_0)}\sin(\rawPhi) } \vee 0.
\end{align*}

Minimization of the resulting function
\begin{align}\label{eq:hed-f2}
\mathscr{F}^{hedge}_2(u; {\tk},S^n_{\tk},W^n_{\tk},\procTh^n_{\tk}) := \sup_{(\mu,\sigma)\in\Theta({\tk},\bar{\mu}^n_{\tk},\procSigma^n_{\tk})}\hat{E}[\hat{V}(\tkp,\mathbf{T}({\tk},x^n_{\tk},u,(\mu,\sigma),\epsilon_{\tkp}))]
\end{align}
with respect to $u\in\cU$ is done by using \texttt{fminbnd} in Matlab with a tolerance threshold of $10^{-6}$, where we double-check whether the minimum is achieved at the boundaries $u \in \{ 0,1\}$.
For the super-hedging loss function $\ell(h) = h^+$, we note that when $W \gg P^{BS}({\tk},S; \procTh)$, the expected loss is numerically zero, so that $V({\tk},S,W)$ is numerically constant and the maximizer $\rawU({\tk},x)$ of \eqref{eq:hed-f2} is ill-defined. This is the scenario where the option can be super-hedged with very high probability so it does not matter what the current hedge is. We found that it helps to detect such cases a priori, setting $\rawU({\tk},x) = 0$ for them.

 \begin{figure}[ht]
   \centering \begin{tabular}{ccc}
     \includegraphics[width=0.32\textwidth]{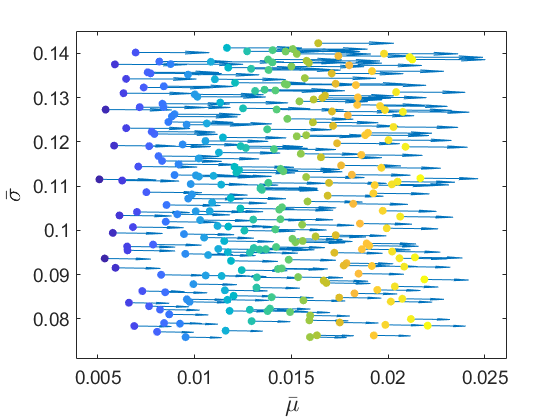} & \includegraphics[width=0.32\textwidth]{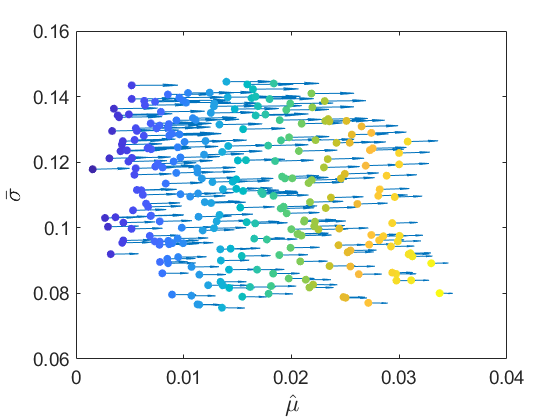} &
     \includegraphics[width=0.32\textwidth]{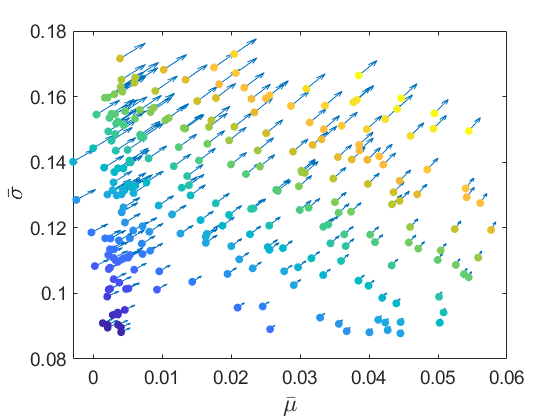} \\
     $t=\Delta t$ & $t = 4 \Delta t$ & $t= 9 \Delta t$ \end{tabular}
     \caption{Estimated angle $\rawPhi$ (quiver direction) and distance $\rawRho$ (quiver length) for the robustified model parameters $\check{\theta}$ as a function of current beliefs $(\procMu,\procSigma)$ for the hedging problem with $S_{t}=100$, $W_{t}=P^{BS}({\tk},100)$. The quivers are color coded in terms of $\modelU({t},\procTh_{t})$ (purple is  $\modelU({t},\procTh)=0$, yellow is $\tilde{u}({t},\procTh)=1$) that is based on a GP surrogate.  \label{fig:hedge-arrows}}
 \end{figure}

\begin{figure}[ht]
   \centering
     \includegraphics[width=0.45\textwidth]{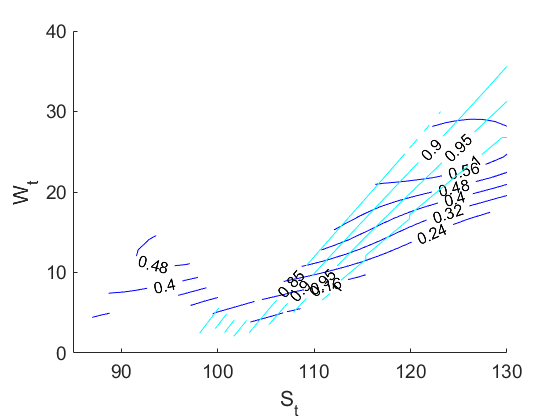}
     \caption{Contour plot of the adaptive robust hedging strategy $\modelU(\tk,S,W, \procMu, \procSigma)$ based on the GP surrogate within the region $ \{(S,W): 0.5P^{BS}({\tk},S)\leq W \leq1.5P^{BS}({\tk},S) \}$ with $\bar{\mu}=0.12$, $\bar{\sigma}=0.4$, $t=t_0$ (blue), $t=t_9$ (cyan). \label{fig:hedge-contour}}
 \end{figure}

 Figure~\ref{fig:hedge-arrows} shows the estimated angle $\rawPhi(t_k, \procMu,\procSigma)$ and distance $\rawRho(t_k, \procMu,\procSigma)$ via a quiver plot at three different time instances $t_k = k \Delta t$. A first observation is that $\rawRho$ is not trivial as the adversarial beliefs are not necessarily on the boundary. Secondly, in this problem $\rawPhi$ depends both on the current stock price $S_{\tk}$ and the current wealth $W_{\tk}$. In the plot we consider at-the-money case $S_{\tk} = \mathcal{K} = 100$ and at-the-wealth case $W_{\tk} = P^{BS}({\tk},S_{\tk})$. We observe that for early $t_k$, the worst-case adversarial $\rawTh$ corresponds to maximizing the stock return $\mu$, $\rawPhi({t_1}, \cdot) \approx 0$, so as to maximize the Call value. Asset volatility begins to play a more important role as expiration approaches and the Vega/Gamma of the Call increase. As a result, $\rawPhi$ rotates counterclockwise and the worst-case $\rawTh$ now involves making \emph{both} $\mu$ and $\sigma$ larger. Third, the Figure indicates the impact of ${\tk}$ on $\modelU({\tk}, \cdot)$: the distribution of $\modelU$ follows the rotation of $\rawPhi$. Time-dependency of $\modelU$ is also illustrated in Figure~\ref{fig:hedge-contour}. The optimal strategy is increasing in $W$ when $W$ is relatively small, and is decreasing when $W$ is large for $t=t_0$. On the other hand, the hedging strategy in general increases with respect to $S$ at $t=t_9$. These show the competition between the two factors $\Phi(S)$ and $W$ at different time steps. At early stages $W$ plays a more important part: small $W$ causes under-hedge and large $W$ leads to over-hedge. When we approach terminal time, $\Phi(S)$ is the main driver of the hedging error: higher stock price implies higher projected option payoff, hence one would invest a larger proportion of the wealth into the stock to reduce the hedging error. In addition, change in the region of the contour reflects different experimental designs $\cD_{k}$ for different time steps.

\subsection{Comparison to Other Approaches}
As comparators to adaptive robust control, we consider again the myopic adaptive (MA) and static robust (SR) formulations, as well as a naive adaptive Delta (AD) hedging. The latter  buys $\Delta({\tk},S_{\tk}; \procTh_{\tk})$ shares at each step ${\tk}$, plugging-in the latest parameter estimates into the classical Delta-hedging strategy of the Black-Scholes  model. MA control corresponds to first solving the optimization problem treating $\theta$ as a parameter, which yields $\tilde{u}^{MA}({\tk}; \theta)$. The adopted strategy is then $\tilde{u}^{MA}({\tk}; \procTh_{\tk})$ which, like the adaptive Delta, leads to a purely learning-based approach. The SR formulation takes $\Theta_{\tk} \equiv \Theta_0$ and then solves the Bellman equations \eqref{eq:hedge-bellman}.

% according to estimated Black-Scholes $\Delta$ as well as 3 control methods: adaptive robust, adaptive, and static robust. Adaptive control is a parametric optimization method where optimal strategy is backward constructed from solving Bellman equation for any fixed $\theta$. The strategy is then seen as a function parameterized by $\theta$ and $t$. At time $t$, the adaptive optimal strategy is the function value at the point estimator of $\theta$. Therefore, it is a purely learning based non-robust approach. Static robust, on the contrary, do not allow belief updating: for any fixed model parameter $\theta$, process of optimal strategy is obtained by solving Bellman equations. Then, the static robust strategy process is the one corresponding to the worst case model.

 \begin{figure}[ht]
   \centering
     \includegraphics[width=0.325\textwidth]{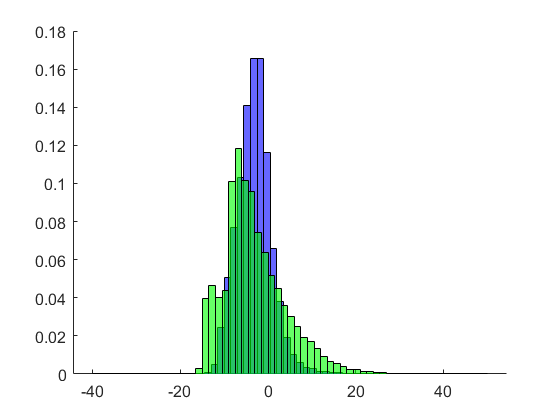}  \includegraphics[width=0.325\textwidth]{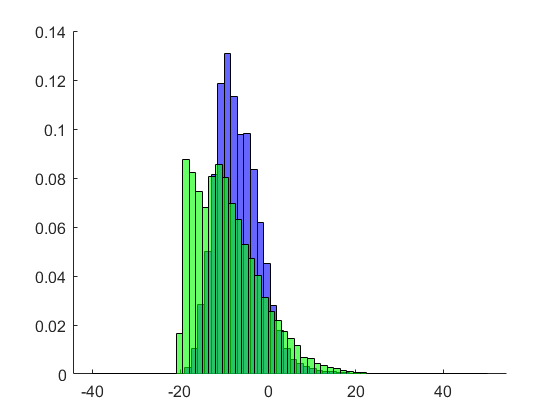}
     \includegraphics[width=0.325\textwidth]{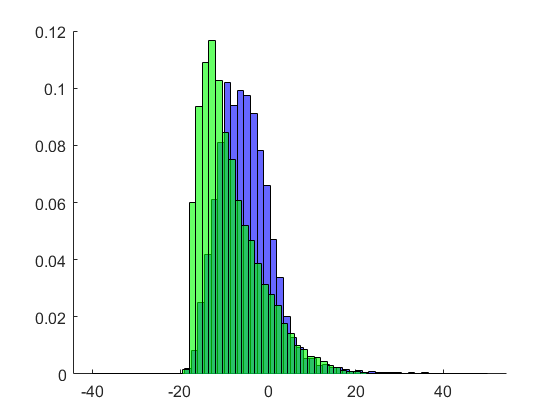}
     \caption{Histograms of hedging error $H=\Phi(S_T)-W_T$ for loss function $\ell(h)=h^++0.75h^-$. Adaptive robust control (blue); adaptive Black-Scholes $\Delta$ (green). Left: $S_{t_0}=90$, $W_0=10$; Center: $S_{t_0}=100$, $W_{t_0}=20$; Right: $S_{t_0}=110$, $W_{t_0}=25$. \label{fig:hist-adap}}
 \end{figure}

To make the comparison, we compute the distribution of the respective terminal hedging errors $H = \Phi(S_T)-W_T$ on $N'=50000$ out-of-sample forward paths, using $\ell(h) =h^++0.75h^-$.
Figure~\ref{fig:hist-adap} shows the empirical histogram of $H$ from the forward simulations.
Better hedging corresponds to lower average loss $\E[ \ell(H)]$ and in particular should concentrate errors closer to zero, so that $\text{mean}(H) \simeq 0$ and $\text{Var}(H)$ small are preferred.
We see in Figure~\ref{fig:hist-adap} that the hedging error of Adaptive Robust strategy is more concentrated around $H=0$. In the out-of-the-money and in-the-money cases, Adaptive Delta leads to significantly larger tail on the positive side. Hence, AR produces a comparatively better strategy which ``super-hedges'' (consistently more negative $H$) AD. This is not very surprising as the AD recipe is ad hoc and is not targeting the loss function.

  \begin{figure}[ht]
   \centering
     \includegraphics[width=0.325\textwidth]{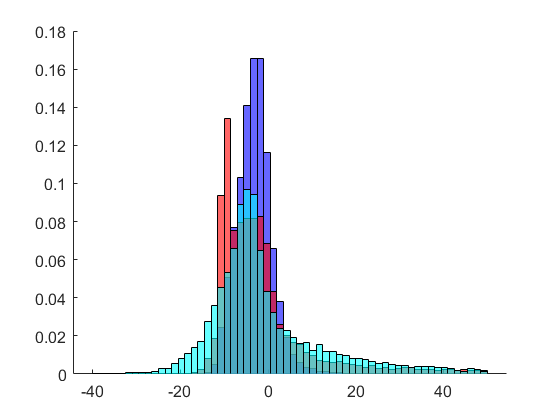}  \includegraphics[width=0.325\textwidth]{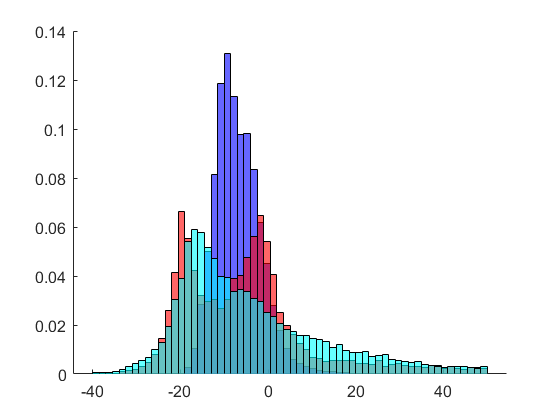}
     \includegraphics[width=0.325\textwidth]{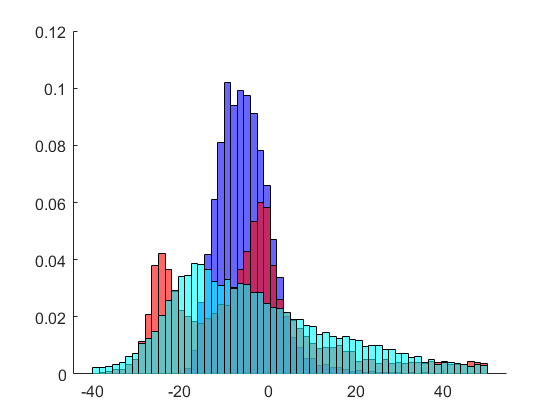}
     \caption{Distribution of $H=\Phi(S_T)-W_T$: adaptive robust (blue); static robust (red);  myopic adaptive (cyan) for loss function $\ell(h)=h^++0.75h^-$. Left: out-of-the-money $S_{t_0}=90$, $W_{t_0}=10$; Center: at-the-money $S_{t_0}=100$, $W_{t_0}=20$; Right: in-the-money $(S_{t_0}, W_{t_0})= (110,25)$. \label{fig:hist-robu}}
 \end{figure}

\begin{table}[tb]
\centering
\renewcommand{\arraystretch}{1.3}
\begin{tabular}{c|rrr||rrr||rrr}
%\cline{2-10} \cline{2-10}
 \multicolumn{1}{c}{ }& \multicolumn{3}{||c||}{$\lambda = 0$}  & \multicolumn{3}{c||}{$\lambda = 0.5$} & \multicolumn{3}{c}{$\lambda = 0.75$}  \\
  %\cline{2-10}
  \multicolumn{1}{c}{ }& \multicolumn{1}{||c}{AR} & SR & MA & AR & SR & MA & AR & SR & MA  \\
  \hline
   mean($H$) & \multicolumn{1}{||c}{-6.09} & -8.82 & -2.79 & -7.11 & -3.34 & -0.60 &  -6.81 & -1.01 & -1.26 \\
 std($H$) & \multicolumn{1}{||c}{16.06} & 5.07 & 21.62 & 5.24 & 19.62 & 25.05 & 5.26 & 24.09 & 22.95  \\
 $q_{0.95}(H)$ & \multicolumn{1}{||c}{21.76} & -1.04 & 37.57 & 1.72 & 33.19 & 49.81 & 1.86 & 52.07 & 43.06  \\
 $V_0$ & \multicolumn{1}{||c}{3.52} & 0.08 & 6.50 & 4.07 & 8.44 & 13.83 & 5.70 & 13.59 & 14.31  \\
 \hline
\end{tabular}
\bigskip
\caption{Mean, standard deviation, and 95\%-quantile of the hedging error $H=\Phi(S_T)-W_T$, as well as the mean loss $V({t_0}, S_{t_0}, W_{t_0})$ with $S_{t_0}=100$, $W_{t_0}=20$. We compare the Adaptive Robust (AR), Static Robust (SR) and Myopic Adaptive (MA) methods and three different loss functions $\ell(h)$ in \eqref{eq:riskmeasure2}.}
\label{table:comparisons}
\end{table}

Figure~\ref{fig:hist-robu} shows that the AR approach frequently performs better than SR and MA strategies. When taking a two-sided loss function ($\lambda=0.75$), the hedging error of AR concentrates more around zero and has apparent thinner tails on both sides compared to the other two. Due to lack of learning, SR strategies are more ``conservative'' and over-invest in the underlying asset (see Figure~\ref{fig:hedge-path} below) which increases the variance of $\Phi(S_T)-W_T$. Myopic adaptive control is over-optimistic and tends to lead to very high variance, already noted in~\cite{bielecki2017adaptive}. %The above observations are verified by numerical computations with different loss functions under these 3 methods in Table~\ref{table:comparisons}.
Table~\ref{table:comparisons} shows the impact of $\lambda$ on the performances of these three approaches. SR does better than the other two when $\lambda=0$. Note that this corresponds to the super-hedging objetive, hence it is not surprising that the conservative approach wins as it over-invests in the risky asset. For $\lambda>0$ where positive hedging error is also penalised, adaptive robust formulation outperforms the other two. Choosing the optimal strategy in these cases are more delicate, and a method which combines learning and robustness handles the situation better.

\subsection{Comparison of Optimal Hedges}

\begin{figure}[ht]
  \centering
    \includegraphics[width=0.45\textwidth]{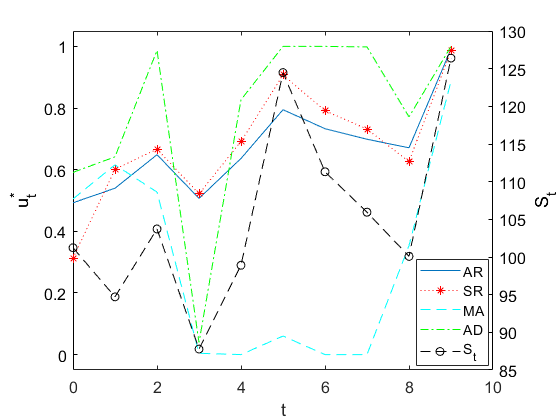}  \includegraphics[width=0.45\textwidth]{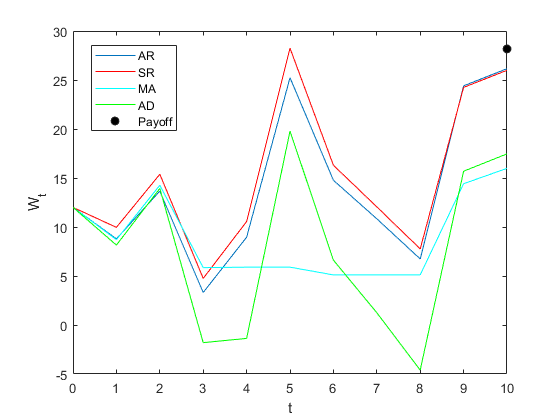}
    \caption{\emph{Left panel:} path of adaptive robust hedge ${\tk} \mapsto \tilde{u}^{AR}({\tk},x_{\tk})$ (blue) for a European Call option in comparison to static robust $\tilde{u}^{SR}({\tk},x_{\tk})$ (red), myopic adaptive $\tilde{u}^{MA}({\tk},x_{\tk})$ (cyan), and adaptive Delta  $\tilde{u}^{AD}({\tk},S_{\tk}; \procTh_{\tk})$ (green). The respective path of $S_{\tk}$ is shown in black (recall strike of $\mathcal{K}=100$). \emph{Right:} paths of the respective wealth processes $W_{\tk}$ color-coded accordingly. The difference between plotted terminal wealth $W_T$ and $\Phi(S_T) \simeq 29$ (black dot, option expires in-the-money) is the respective hedging error $H$.
\label{fig:hedge-path}}
\end{figure}

Compared to hedging strategies generated by other  approaches, the adaptive robust strategy follows the movement of the underlying $S_t$ more closely. Such behavior can be understood as a better learning of the model than other methods. Figure~\ref{fig:hedge-path} displays a sample trajectory of $S_t$ and corresponding $\tilde{u}({\tk}, S_{\tk})$ across the above three approaches.
When stock price is volatile during the time period, robust strategies $\tilde{u}^{AR}({\tk},S_{\tk})$ and $\tilde{u}^{SR}({\tk},S_{\tk})$ are much more stable than the pure learning-based strategies $\tilde{u}^{MA}({\tk},S_{\tk})$ and $\tilde{u}^{AD}({\tk},S_{\tk},\bar{\theta}_{\tk})$.
%When stock price stays deep in the money $S_t \gg K$, $\tilde{u}^{AR}(t,S_t)$ is still sensitive to further changes in $S_t$, even though the Black-Scholes $\Delta(t,S_t) \simeq 1$ will be constant and effectively asset-independent.
It is also worth mentioning that for the static robust approach the assumed drift $\check{\mu}$ corresponding to the worst case model is very high. Namely the SR worst case is that the Call ends in-the-money and a large positive hedging error $H \gg 0$ results. To compensate against this scenario, the static robust strategy $\tilde{u}^{SR}$ over-invests in the risky asset to ensure that $W_T \simeq \Phi(S_T)$ conditional on $S_T \gg \mathcal{K}$, which conversely generates larger risk in other, less adversarial scenarios. This explains why $u^{SR}({\tk},S_{\tk}) > u^{AR}({\tk},S_{\tk})$ in Figure~\ref{fig:hedge-path}.

% \begin{figure}[ht]
%   \centering
%     \includegraphics[width=0.325\textwidth]{RegressionPictures/hedging_struct_o.png}  \includegraphics[width=0.325\textwidth]{RegressionPictures/hedging_struct_a.png}
%     \includegraphics[width=0.325\textwidth]{RegressionPictures/hedging_struct_i.png}
%     \caption{Adaptive robust hedge $u^*$ as a function of $\hat{\mu}, \hat{\sigma}, S$ for a European Call options. We show three cuts for three different levels of $S_t$. The blue horizontal plane indicates the Black-Scholes $\Delta(t,S_t)$. \blu{what is $W_t$ and $t$? \label{fig:hedge-surf}}
%     }
% \end{figure}

\section{Discussion}\label{sec:enhance}

\subsection{Enhancements}
A central feature of our machine learning approach is the rich set of opportunities to enhance the computations. These concern the ultimate aim of (i) faster running time; (ii) more stable and accurate estimates of $\hat{V}$ and $\modelU$. Practically, we would like to be able to specify a small simulation budget $N_k$ and then obtain good results quickly. We do note that there is a third dimension of actual running time---if the statistical surrogate model is too complicated then its overhead can be high even when $N_k$ is low. For example, we can propose heuristics for sequential, adaptive designs $\cD_{k}$ that would maximize the learning ability of $\modelU$ as a function of $N_k$, but actually slow down the algorithm (hence we deem them low-priority extensions).

%, i.e.~lower MC errors. Practically, we measure aspect speed through the simulation budget $N$ which we try to make as small as possible.

\textbf{Faster optimization:} Generating a single training pair $(x_{\tk}, v_{\tk})$ in \eqref{eq:yn} requires solving a nested optimization problem (maximizing over $u$, minimizing over $(\mu,\sigma)$), which is the most expensive part of the algorithm. Since this is done inside all the loops (over $k$ and over $n$), making this step efficient is central to fast performance. In our existing implementation  we utilize off-the-shelf gradient-free optimizers (\texttt{fminbnd} and \texttt{fminsearch} in \texttt{Matlab}) which are agnostic to the structure of $\hat{E}$ and $\hat{V}(\tkp,\cdot)$ and operate completely generically. One way to speed up isthrough warm-starting the optimizers at the location $\rawU(\tkp, x^n_{\tk})$, i.e.~assuming that the new optimizer is the same as the one from the previous (later) time-step. For $\mathscr{F}_2$ in \eqref{eq:f2}-\eqref{eq:hed-f2} we may warm-start the optimizer at $\varphi = \pi$. A further enhancement is to rely on gradient-descent optimizers which would require to estimate $\nabla \mathscr{F}_2(u)$. Within the GP setup we may analytically differentiate the GP kernel to learn $\nabla_x \hat{V}$ which in turn could be converted into $\nabla_u \mathscr{F}_2$.

%Show a table with number of \texttt{predict} calls using different tricks: fminbnd (which has no warm starts), fminsearch (which does) -- with different ways to make guesses.

\textbf{GP Hyperparameters:} A significant portion of the algorithm running time is spent on fitting the GP surrogates, i.e.~inferring (typically by maximum likelihood maximization) the respective hyperparameters $\vartheta_{\tk}$. This is a hard nonlinear optimization problem with total time complexity $\mathcal{O}(K \cdot N^3)$, with the dependence on $K$ due to the need to fit a surrogate at each time-step $\tk$. Generally, the hyperparameters $\hat{\vartheta}_{\tk}$ are stable over time, so one could use this for another warm-start in the MLE optimizer, or even directly set $\hat{\vartheta}_{\tk} = \hat{\vartheta}_{\tkp}$, freezing hyperparameters across (some) time-steps. Because the GP surrogate is data-driven, the final prediction depends mostly on the training data and less on the precise value of $\vartheta_{\tk}$.

In terms of accuracy, an important piece of building the surrogate are its extrapolation properties. Recall that for extrapolation the GP prediction reverts to its prior $m_*(x) \to m_0(x)$ outside of the input domain. For fitting the control surrogate $\modelU$ in both case studies we take $m_0(x) \equiv 0$, but other choices would be appropriate depending on the problem/surrogate and in our experience can play a nontrivial role in ultimate solution quality. (However, note that $m_0(\cdot)$ is only relevant for extrapolation---within $\cD$ the prediction is driven by the training data.)

\subsection{True Monte Carlo}
The use of numerical quadrature to integrate out $\hat{V}(\tkp,\cdot)$ is probabilistically biased in the sense that $\hat{E}^{Quad}[ \hat{V}(\tkp,\cdot) ](x) \neq \E[ \hat{V}(\tkp,\cdot)](x)$ so that a consistent, deterministic error is introduced during this sub-step. A well-known alternative is to employ Monte Carlo integration, namely replacing the integral with a \emph{random} sum
\begin{align}
  \int \hat{V}(\tkp,\bT({\tk},x,u,\theta,z)) f_\eps(z) dz \approx \frac{1}{\mathfrak{I}} \sum_{i=1}^{\mathfrak{I}} \hat{V}(\tkp, \bT({\tk},x,u,\theta,Z_i) ) =: \hat{V}^{MC}({\tk}),
\end{align}
where $Z_i\sim \cN(0,1), i = 1, \ldots, \mathfrak{I}$ are now i.i.d.~Gaussian samples, furthermore independent across different $x^n$'s. Relative to quadrature, the MC approximation   is unbiased, but noisy. We find that the latter noise makes a Monte Carlo-based scheme much less accurate.

%Observation noise is intrinsic to GP surrogates which already estimate $\sigma_\eps$.
Specifically,  for the optimal hedging problem, we computed $\tilde{u}^{\text{MC}}({\tk},\cdot)$ and $\tilde{u}^{\text{Quad}}({\tk},\cdot)$ using several different quantization/inner-simulation levels $\mathfrak{I}$ and relying on identical design
 $(S^{1:N},\bar{\mu}^{1:N},\bar{\sigma}^{1:N},W^{1:N})$. We find that to achieve comparable accuracy one needs 3-5 times more MC points compared to quadrature points. Given that the computational complexity of both approaches is the same, it seems much preferable to employ quantization.

\subsection{Related Control Formulations}
As already demonstrated in the one-period example of Section \ref{sec:one-period}, our numerical tools can be straightforwardly applied to \emph{parametric} control problems. In the latter formulation, rather than seeking to optimize/robustify against a collection of potential $\theta$'s, we simply wish to understand the dependence of the value function or feedback control on $\theta$. The proposed Gaussian process surrogates is a natural tool for that purpose, allowing the modeler to solve the underlying control problem at several instances $\procTh^n$, $n=1,\ldots, N$ and then interpolate to obtain $V(\cdot; \theta)$ and $u^*(\cdot; \theta)$ at arbitrary $\theta \in \Theta$. In particular, such parametric control is relevant for the myopic adaptive problem where one pre-computes $V(\cdot; \theta)$ and then plugs-in the latest belief $\procTh_t$.

Secondly, our scheme also immediately nests the dynamic adaptive problem where the adversarial set $\Theta(\procTh) = \{\procTh\}$ is a singleton, so the inner optimization disappears. Indeed, one can simply take $\kappa = 0$ to collapse the radius of $\Theta$ to zero, and then run exactly the same code as for the adaptive robust setting. Third, the scheme can be adapted to handle Bayesian formulations that replace $\inf_\theta$ with an integral $\int_\theta (\cdot) \nu_{\tk}(d\theta)$ against a distribution $\nu_{\tk}(d\theta)$. Learning becomes encoded as updating the posterior $\nu_{\tk}(\cdot)$ of $\theta$ over time. In the practically tractable case, there are finite-dimensional sufficient statistics (such as the posterior mean) $\procTh_{\tk}$ for $\cL(\theta^*  | \cF_{\tk})$ so the integral over $\nu_{\tk}(d\theta)$ can be reduced to a numerical quadrature over the parametrized posterior distribution. For example, under drift uncertainty and assuming the setting of trying to learn an unknown fixed $\mu^*$, $\cL(\mu^* | \cF_T)$ is Gaussian with a certain formula for mean $\procMu_{\tk}$ and posterior variance $V_t$ and integrating against $\nu_{\tk}(d\theta)$ becomes a one-dimensional Gaussian integral.

\begin{remark}
Our algorithm can also be interpreted as a type of Control Randomization. In CR, the solution pipeline is of the form depicted below, where $\vec{u}^{(0)}$ is possibly randomized.
%\tikzstyle{decision} = [diamond, draw, fill=blue!20,
%    text width=4.5em, text badly centered, node distance=3cm, inner sep=0pt]
\tikzstyle{block} = [rectangle, draw, fill=blue!20,
    text width=7em, text centered, rounded corners, minimum height=4em,node distance=4.5cm]
\tikzstyle{line} = [draw, -latex']
\tikzstyle{cloud} = [draw, ellipse,fill=red!20, node distance=3cm,
    minimum height=2em,text centered, text width=5em]

\bigskip

\begin{tikzpicture}[node distance = 3cm, auto,scale=0.7,every node/.style={transform shape}]
    % Place nodes
    \node [block] (init) {Propose a control strategy $\vec{u}^{(0)}$ and measure $\Q^{0}$};
    %\node [cloud, above of=init] (rand) {possibly randomized};
    \node [block, right of=init] (sim1) {Simulate  $N$ paths under $\Q^{0}$};
    \node [block, right of=sim1] (design) {Create designs $x^{1:N}_{\tk}$};
    \node [block, right of=design] (dp) {Solve for $\rawU(\tk, x^{1:N}_{\tk})$ via backward induction};
    \node [block, below of=dp, node distance=3cm] (out) {Feedback control map $\modelU(\tk, \cdot)$};
    \node [block, right of=dp] (sim2) {Forward Monte Carlo using $\modelU$ and $\Q^{test}$};
    \node [block, below of=sim2, node distance=3cm] (test) {Estimated value function $\hat{V}(0,x_0)$};
    % Draw edges
    \path [line] (init) -- (sim1);
    \path [line] (sim1) -- (design);
    \path [line] (design) -- (dp);
    \path [line] (dp) --(sim2);
    \path [line] (sim2) -- (test);
    \path [line,dashed] (dp) -- (out);
    %\path [line,dashed] (rand) -| (sim1);

\end{tikzpicture}
\bigskip

With this interpretation, any experimental design can be viewed as an implicit recipe for $\Q^0$, namely specifying the respective marginal distribution of $X_{\tk}$. This perspective emphasizes the fact that the quality of the design is linked to the quality of the ansatz for $u^{(0)}$. It also implies that the training procedure can be embedded into top-level iterations, whereby better and better estimates of $u^{*}$ are fed back as new guesses $u^{(r)}, r=1,\ldots$ to generate more and more targeted simulation designs $x^{1:N,(r)}_{\tk}$.
\end{remark}

\subsection{Conclusion}\label{sec:conclusion}

We have developed a machine learning-based algorithm for adaptive robust control. Our motivation comes from the multiple desirable features of adaptive control and aims to mitigate the associated computational challenges, making the adaptive robust framework numerically feasible. The resulting case studies provide new insights on the interplay between learning, controlling for model uncertainty and risk aversion in the context of the classical Merton problem and the loss-based hedging problem.

The key innovation in our methodology is to build multiple surrogates for different pieces of the Bellman recursion, in particular not just for the value function but also for the feedback control map and worst-case parameters. To do so we propose utilization of Gaussian Process surrogates which check off multiple important computational considerations. The developed algorithm is certainly of independent interest. Further investigations of non-financial contexts and of other related control formulations are deferred to future research.

\bibliographystyle{siam}
\bibliography{optInvInterp}

\end{document}